\documentclass[preprint,12pt,3p]{elsarticle}

\usepackage{lipsum}
\makeatletter
\def\ps@pprintTitle{%
 \let\@oddhead\@empty
 \let\@evenhead\@empty
 \def\@oddfoot{}%
 \let\@evenfoot\@oddfoot}
\makeatother

\usepackage{graphicx}
\usepackage{hyperref}   
\usepackage{lineno}	
\usepackage{amssymb} 
\usepackage{scalerel, amsmath} 
\usepackage{bm}
\usepackage{dsfont}	
\usepackage{navigator}	
\usepackage{changepage}	
\usepackage[pdftex]{color}
\usepackage{caption}
\usepackage{subcaption}
\usepackage{lipsum}
\usepackage{url}
\usepackage{afterpage}
\usepackage[normalem]{ulem}
\usepackage[title]{appendix}

\usepackage{pgfplots}
\pgfplotsset{compat=newest}
\usetikzlibrary{plotmarks}
\usetikzlibrary{arrows.meta}
\usepgfplotslibrary{patchplots}
\usepackage{grffile}

\newcommand{\vxi}{\boldsymbol{\xi}}

\newcommand{\beq}{\begin{linenomath*} \begin{equation}}
\newcommand{\eeq}{\end{equation} \end{linenomath*}}
\newcommand{\bseq}{\begin{linenomath*} \begin{subequations} }
\newcommand{\eseq}{\end{subequations} \end{linenomath*}}

\newcommand{\bal}{\begin{align}}
\newcommand{\eal}{\end{align}}

\newcommand{\bfi}{\begin{figure}}
\newcommand{\efi}{\end{figure}} 

\newcommand*\colvec[3][]{
    \begin{pmatrix}\ifx\relax#1\relax\else#1\\\fi#2\\#3\end{pmatrix}
}

\graphicspath{{./figures/}}

\newcommand{\code}[1]{\texttt{#1}}
\usepackage{booktabs}

\journal{Elsevier}

\begin{document}

\title{Solid shell prism elements based on hierarchical, heterogeneous, and anisotropic shape functions}

\author[1]{\L ukasz Kaczmarczyk*}

\author[1]{Hoang Nguyen} 

\author[2]{Zahur Ullah}

\author[1]{Mebratu Wakeni}

\author[1]{Chris Pearce}

\address[1]{Glasgow Computational Engineering Centre, James Watt
School of Engineering, University of Glasgow,
Glasgow, G12 8QQ, UK}

\address[2]{Advanced Composites Research Group, School of 
Mechanical and Aerospace Engineering, Queen's University 
Belfast, Belfast, BT9 5AH, UK}

\begin{abstract}
The formulation of a new prism finite element is presented for
the nonlinear analysis of solid shells subject to large strains and large displacements. 
The element is based on hierarchical, heterogeneous, and anisotropic
shape functions. As with other solid shell formulations, only displacement 
degrees of freedom are required to describe the shell kinematics and general 
three-dimensional material laws can be adopted.
However, the novelty of this formulation is the ability to capture complex shell behaviour 
and avoid locking phenomena, without the need to use reduced integration or 
adopt additional natural strain or enhanced strain fields. 
Thus, this element is ideally suited for geometrically and physically nonlinear problems.
This is achieved by constructing independent approximation shape functions on both the prism element's triangular faces
and through the thickness, where the latter is associated 
with a local coordinate system that convects with deformation of the shell. 
The element is extremely efficient, with the hierarchical property lending itself to an efficient and highly scalable multigrid solver,
and the heterogeneity property enables local p-adaptivity. The paper demonstrates performance of the element for a number of linear and geometrically nonlinear problems, benchmarked against well established problems in the literature.
The formulation has been implemented in the \code{MoFEM} library \cite{kaczmarczyk_mofem_2020}. Both the code and the data
for the numerical examples are open-source \cite{mofem_solid_shell_code_and_data}.
\end{abstract}

\begin{keyword}
  hierarchical shape functions \sep 
  prism element \sep 
  solid shell finite element \sep
  large deformations \sep
  multi-grid
\end{keyword}



\maketitle

\setlength\textwidth{177.8mm}

\section{Introduction}
\label{sec:introduction}

Shell structures continue to generate scientific interest and represent an ongoing challenge in terms of computational mechanics.  Shells are typically thin, curved and lightweight structures, and are technically important across a range of applications, including civil and structural engineering, aerospace, automotive, industrial processing and marine. 
There is also a growing interest in the mechanics of shells at the interface of engineering and the life-sciences, where they occur as natural structures comprised of complex materials with unusual nonlinear behaviour.

Degenerated shell theory has been a successful approach for modelling shells for several decades and has worked well for many cases. However, there are a number of restrictions and scenarios where it is less effective. The kinematic assumptions and the rotational degrees of freedom can lead to difficulties. Although in-plane stresses and transverse shear stresses are considered, the stretching effect in the through-thickness direction is not taken into account. Consequently, it is necessary to reduce the general three-dimensional material law to the plane stress condition. When describing the boundary conditions for geometrically nonlinear problems, it is necessary to undertake a complex update of the rotations. Problems also arise when shell elements are used in combination with solid elements. Therefore, it is
desirable to have an alternative approach that is able to overcome these shortcomings.

Solid shell theory
\cite{parisch_continuum-based_1995, hauptmann_systematic_1998} assumes the
element has only displacement DOFs. As a result, these elements are capable of connecting naturally with
continuum elements and overcome the disadvantages related to boundary conditions that are inherent in degenerated shell elements.
There is a well-established
body of work on the linear and nonlinear analyses of shells using solid shell theory. 
Hosseini et al. \cite{hosseini_isogeometric_2013}
employed it for nonlinear analysis where the isogeometric
formulation was used to construct the shell's
mid-surface and linear Lagrange function for the approximation through the
thickness. Leonetti et al. \cite{leonetti_efficient_2018} also worked on the
isogeometric implementation of a solid shell element for geometrically nonlinear
problem where a generalised constitutive matrix is employed. Solid shell elements
also find applications in modelling laminated composites and delamination
phenomenon \cite{hashagen_numerical_2000,remmers_solid-like_2003,
vu-quoc_optimal_static_2003,vu-quoc_optimal_dynamic_2003}. 

Despite the advantages of solid shell formulations, 
locking phenomena, such as shear locking, membrane locking, and especially volumetric and
thickness locking, are common problems
\cite{sze_three-dimensional_2002, magisano_advantages_2017}. There has been
considerable efforts to circumvent these by means 
of assumed natural strain, enhanced assumed strain, and mixed
stress-displacement formulations.
For example, in order to overcome thickness and volumentric locking, a hybrid-stress formulation was employed~\cite{sze_hybrid_part1_2000}. 
Reese et al. \cite{reese_new_2000} proposed a
brick element formulation with the combination of reduced integration and
stabilisation concept for large deformation problems based on the enhanced strain
method. Later, Reese~\cite{reese_large_2007} also proposed an eight-node solid
element with reduced integration based on hourglass stabilisation for static
problems.  
Schwarze and Reese \cite{schwarze_reduced_2009,schwarze_reduced_2011}
developed a reduced integration eight-node solid-shell finite element, requiring only one 
integration point within the shell plane and at least two integration points through the thickness. 
The formulation has the advantage of being able to choose an arbitrary number of integration 
points through the shell thickness and only one enhanced degree-of-freedom is needed to avoid 
volumetric and thickness locking. The assumed natural strain concept is adopted to avoid 
shear and curvature thickness locking. 


In the present study, a new solid shell prism element is proposed. 
Independent approximation shape functions are constructed on the triangular faces
and through the thickness of the prism element, where the latter is associated 
with a local coordinate system that convects with deformation of the shell. 
This anisotropic property of the shape functions avoids locking, without the need to use an assumed
strain field or enhanced assumed strain field, or to use reduced integration. 
Moreover, the prism element is constructed using hierarchical and heterogeneous shape functions, 
enabling local $p$-refinement and development of a multigrid solver \cite{mitchell_hp-multigrid_2010}, that 
means the element is also computationally efficient.

The outline of this paper is as follows. Section \ref{sec:prism_element}
presents the framework to construct hierarchical approximations for
triangular prism elements and Section \ref{sec:solid_shell} discusses the
kinematics and linearisation following the solid shell theory. The geometry and
displacement approximations and discretisations are given in Section
\ref{sec:aproximations}. Section \ref{sec:example} evaluates the convergence 
properties, as well as presenting some benchmark examples of shells with large
deformations. Concluding remarks are presented in
Section \ref{sec:conclusion}.

\section{Construction of hierarchical shape functions for prism elements}
\label{sec:prism_element}

\begin{figure}[!htbp]
	\centering
	\includegraphics[scale=0.5]{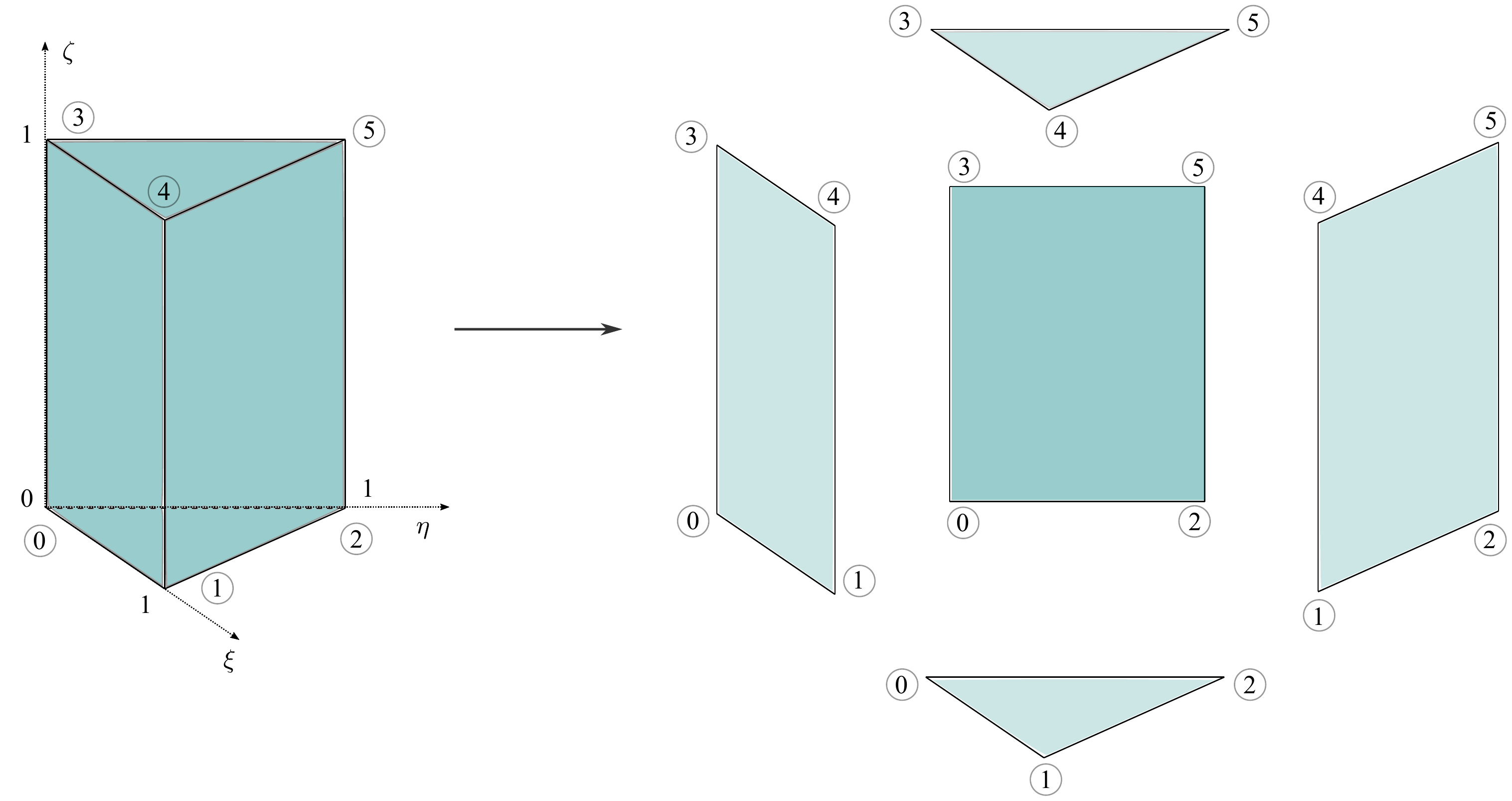}
	\caption{Geometry of prism element and exploded view with node numbering in circles.}
	\label{fig:prism_geometry}
\end{figure}

This section presents the construction of the shape function for the proposed solid shell prism element, shown in Fig. \ref{fig:prism_geometry}. This
approximation is hierarchical, heterogeneous, and anisotropic. The
hierarchy of the approximation is built on different geometrical entities, i.e.
vertices, edges, faces, and volume. 
Meanwhile, an arbitrary order of approximation can be established independently 
for each geometrical entity, including each prism, enabling heterogeneous approximations across a finite element mesh. 
Moreover, approximations functions on the triangular faces are defined
independently from those of the quadrilateral faces.

The starting point for construction of the element is
the \textit{primary basis} from a set of one-dimensional polynomials $\left\{ {{\psi_\ell}:\ell=
0,1,\dots} \right\}$ that constitute a hierarchic basis on a reference interval
$\left[-1,1\right]$ \cite{ainsworth_hierarchic_2003}. The function $\psi_\ell$
is usually in the form of Legendre polynomial $L_\ell$ with degree $\ell$ which
are defined as follows 
\beq\label{eq:Lagrange}
{L_0}\left( s \right) = 1,  \quad {L_1}\left( s \right) = s,  \quad {L_{\ell + 1}}\left( s \right) = \frac{{2\ell + 
1}}{{\ell + 1}}s{L_\ell}\left( s \right) - \frac{\ell}{{\ell + 1}}{L_{\ell -
1}}\left( s \right), \quad s \in \left[ -1, 1\right], \quad \ell = 
1,2,\dots
\eeq
It is worth noting that a more general form of Legendre polynomials presented here, such as Gegenbauer polynomials, can be found in the literature \cite{gradshteyn_table_1994,ainsworth_hierarchic_2003}.

A prism element consists of coordinate points $(\xi, \eta, \zeta)$ such that 
$0\leq \xi, \eta, \zeta \leq 1$ and $\xi + \eta \leq 1$. In other words, an element is the 
Cartesian product of a triangle element with coordinates $(\xi, \eta)$ and a one-dimensional segment element
identified by coordinate $\zeta$. Thus, barycentric 
coordinates $\lambda_j$ ($j = 0, 1, 2$) of a triangle and affine 
coordinates $\mu_i$ ($i = 0, 1$) of a one-dimensional
segment are defined as follows
\bseq
\begin{align}
\lambda_0 &= 1 - \xi - \eta, \quad \lambda_1 = \xi, \quad \lambda_2 = \eta, \\
\mu_0 &= 1 - \zeta, \quad \mu_1 = \zeta, 
\end{align}
\eseq
and their gradient with respect to the Cartesian coordinates $(\xi, \eta,
\zeta)$ given by
\bseq
\begin{align}
  \nabla\lambda_0 &= \colvec[-1]{-1}{0}, \quad \nabla\lambda_1 = \colvec[1]{0}{0}, \quad \nabla\lambda_2 = \colvec[0]{1}{0}, \\
  \nabla\mu_0 &= \colvec[0]{0}{-1}, \quad \nabla\mu_1 = \colvec[0]{0}{1}. 
\end{align}
\eseq

\subsection{Vertex shape functions}
The shape function for a given vertex $\mathrm{v}$ (v$ = 0,\dots, 5$), see Fig. \ref{fig:prism_geometry}, is a bi-linear 
polynomial that takes the value of one at the vertex to which it is associated and zero at the other vertices. 
This, together with the fact that it is linear along edges, is 
sufficient to satisfy the necessary $C^0$-continuity of the global vertex functions on adjacent prism elements 
that share the vertex. Thus the shape 
function $\phi^{\mathrm{v}}$ 
is expressed as the product of the triangle's barycentric coordinates and the segment's affine coordinates as follows
\beq
\phi^\mathrm{v} = {\lambda _i}\mu_j,
\eeq
where the relationship between the indices $\mathrm{v},i$ and $j$ is given in Table \ref{tab:vertex_bases}

\begin{table}[htbp]
  \centering
  \caption{Relationship between the indices of the vertex shape functions}
    \begin{tabular}{lll}
    \hline
    $ \mathrm{v} \quad $     & $i \quad $       &  $j \quad $  \\
    \hline
    0     & 0     & 0 \\
    1     & 1     & 0 \\
    2     & 2     & 0 \\
    3     & 0     & 1 \\
    4     & 1     & 1 \\
    5     & 2     & 1 \\
    \hline
    \end{tabular}
  \label{tab:vertex_bases}
\end{table}

\subsection{Edge shape functions}
There are two types of edges to be considered: triangle edges and through-thickness edges. 
The Lagrange polynomials given by Eq. \eqref{eq:Lagrange} are used to define the shape functions for both of these.

Consider first the triangle edges of the prism's top and bottom triangular faces. 
The Lagrange polynomials are mapped to an edge with the interval $\left[-1,1\right]$ using the difference of the barycentric coordinates as $L_{\ell}(\lambda_j - \lambda_k)$.
Care must be taken to use the global (rather than local) number ordering
of the edge, since edges that are shared by adjacent prism elements may have different local numbering. 
For example, if the local and global number ordering of Edge 0-1 are the same, then the mapping $\lambda_1 - \lambda_0$ is used, otherwise the difference is reversed. 

The shape function for any edge is non-zero on the edge to which it is associated and zero on all
other edges. For Edge 0-1, for example, the function $\mu_0\lambda_0\lambda_1$ satisfies this property.
Hence,
assuming the same local and global numbering, the shape functions associated with
edge $k$-$j$ are
\beq
\phi^{\mathrm{te}}_\ell = \beta^{\mathrm{te}}_{ijk} L_{\ell}(\lambda_j - \lambda_k),
\eeq
where $\beta^{\mathrm{te}}_{ijk} = \mu_i\lambda_j\lambda_k$ and the relationship between
indices is given in Table \ref{tab:triangle_edge_bases}. $\ell
= 0, \dots, p-2$, where $p$ is overall polynomial order of the prism. The
gradients of $\phi^{\mathrm{te}}_\ell$ with respect to $(\xi, \eta, \zeta)$ are
given as
\beq
\begin{aligned}
\nabla\phi^{\mathrm{te}}_\ell = & \nabla\mu_i\,\lambda_j\lambda_k L_{\ell}(\lambda_j - \lambda_k) + \\
                               & \mu_i\nabla\lambda_j\,\lambda_k L_{\ell}(\lambda_j - \lambda_k) + \\
                               & \mu_i\lambda_j\nabla\lambda_k\, L_{\ell}(\lambda_j - \lambda_k) + \\
                               & \mu_i\nabla\lambda_j\,\lambda_k L'_{\ell}(\lambda_j - \lambda_k)\,[\nabla\lambda_j - \nabla\lambda_k].
\end{aligned}
\eeq

\begin{table}[htbp]
  \centering
  \caption{Relationship between indices for triangle edge shape functions}
    \begin{tabular}{llll}
    \hline
    $\mathrm{te} \quad $    & $i \quad $     & $j \quad $     & $k \quad$ \\
    \hline
    0-1   & 0     & 1     & 0 \\
    1-2   & 0     & 2     & 1 \\
    2-0   & 0     & 0     & 2 \\
    3-4   & 1     & 1     & 0 \\
    4-5   & 1     & 2     & 1 \\
    5-3   & 1     & 0     & 2 \\
    \hline
    \end{tabular}
  \label{tab:triangle_edge_bases}
\end{table}

The through-thickness edges of the prism's quadrilateral faces are parameterised by the affine coordinates, such that the shape function is expressed as
\beq
\label{eq:phi-qe}
\phi^{\mathrm{qe}}_\ell = \beta^{\mathrm{qe}}_i L_{\ell}(\mu_1 - \mu_0),
\eeq
where $\ell = 0,\dots,p-2$ and the index $i$ takes values 0, 1, and 2 for Edges
0-3, 1-4, and 2-5, respectively. The Lagrange polynomials $L_{\ell}$ have been mapped using the affine coordinates and function $\beta^{\mathrm{qe}}_i=\lambda_i\mu_0\mu_1$ is non-zero on the edge under consideration and zero on all other edges. As above, it has been assumed the local and global number ordering are the same.



\subsection{Face shape functions}
The prism element comprises two types of faces: triangular and quadrilateral. Following a similar approach to edges, the shape functions for triangular faces are expressed as
\beq
\label{eq:phi-triangle}
\phi^{\triangle}_{k\ell} = \beta^\triangle_i L_k(\lambda_1 -
\lambda_0) L_\ell(\lambda_2 - \lambda_0),
\eeq
where $\beta^\triangle_i = \mu_i\lambda_0\lambda_1\lambda_2$ is selected to
enforce the condition that $\phi^{\triangle}_{k\ell}$ is non-zero on the face to
which it is associated and zero on all other faces. The index $i$ takes the
value of 0 and 1 for Triangles 0-1-2 and 3-4-5, respectively.

For quadrilateral faces, the shape functions are expressed as
\beq
\phi^{\square}_{k\ell} = \beta^\square_{ij} L_k(\mu_1 - \mu_0)L_\ell(\lambda_j - \lambda_i),
\eeq
where $\beta^\square_{ij} = \mu_0\mu_1\lambda_i\lambda_j$ and the relationship between the indices are given in Table \ref{tab:quad_face_bases}.

\begin{table}[htbp]
  \centering
  \caption{Relationship between indices for quadrilateral face shape functions}
    \begin{tabular}{lll}
    \hline
    $\square \quad$     & $i \quad$     & $j \quad$ \\
    \hline
    0-1-4-3 & 0     & 1 \\
    1-2-5-4 & 1     & 2 \\
    2-0-3-5 & 2     & 0 \\
    \hline
    \end{tabular}
  \label{tab:quad_face_bases}
\end{table}

\subsection{Volume shape functions}
To construct the volume shape functions, the hierarchical Lagrange polynomials adopted for the triangular faces in Eq.~(\ref{eq:phi-triangle})
are expanded through the thickness of the element using the Lagrange polynomials adopted for the quadrilateral edges in Eq.~(\ref{eq:phi-qe}). Thus, 
\beq
\phi^\mathrm{P}_{klm} = \beta^\mathrm{P} L_k(\mu_1 - \mu_0)L_\ell(\lambda_1 - \lambda_0)L_m(\lambda_2 - \lambda_0).
\eeq
where $\beta^\mathrm{P} = \mu_0\mu_1\lambda_0\lambda_1\lambda_2$ enforces the condition that the volume shape function is zero on all faces.

It will be shown in the next section of this paper that displacement DOFs associated with triangle edges and faces are expressed in the global Cartesian coordinate system, whereas DOFs associated with quadrilateral edges and faces, and prism volumes are expressed in the local, curvilinear, through-thickness coordinate system. This local coordinate system convects with deformation of the shell. The distribution of
DOFs is summarised in Table \ref{tab:dofs_distribution}.

\begin{table}[htbp]
	\centering
	\caption{DOFs distribution in different coordinate systems}
	  \begin{tabular}{lll}
	  \hline
	  DOFs  & \multicolumn{1}{l}{Coordinate system} &  \\
  \cline{2-3}          & \multicolumn{1}{l}{Global Cartesian} & Local curvilinear \\
	  \hline
	  Vertices & \multicolumn{1}{l}{\checkmark} &  \\
	  Edges (triangle) & \multicolumn{1}{l}{\checkmark} &  \\
	  Triangles & \multicolumn{1}{l}{\checkmark} &  \\
	  Edges (quadrilateral) &       & \checkmark \\
	  Quadrilaterals &       & \checkmark \\
	  Volume &       & \checkmark \\
	  \hline
	  \end{tabular}
	  \label{tab:dofs_distribution}
	\end{table}
  
To conclude this section, it is worth noting that the proposed shape functions enable significant flexibility in building approximations across a finite element
mesh. The heterogeneous property allows different polynomial orders
to be applied for different entities at different parts of the mesh. 
The anisotropic property enables different orders to be used in different directions, e.g.
on the top and bottom surfaces or in the through-thickness direction. 
Furthermore, the hierarchic nature of the shape functions
supports effective implementation of a multi-grid solver, which significantly improves
the analysis speed due to its inherent scalability.

\section{Kinematics of solid shells}
\label{sec:solid_shell}

\begin{figure}[!htbp]
	\centering
	\includegraphics[scale=.5]{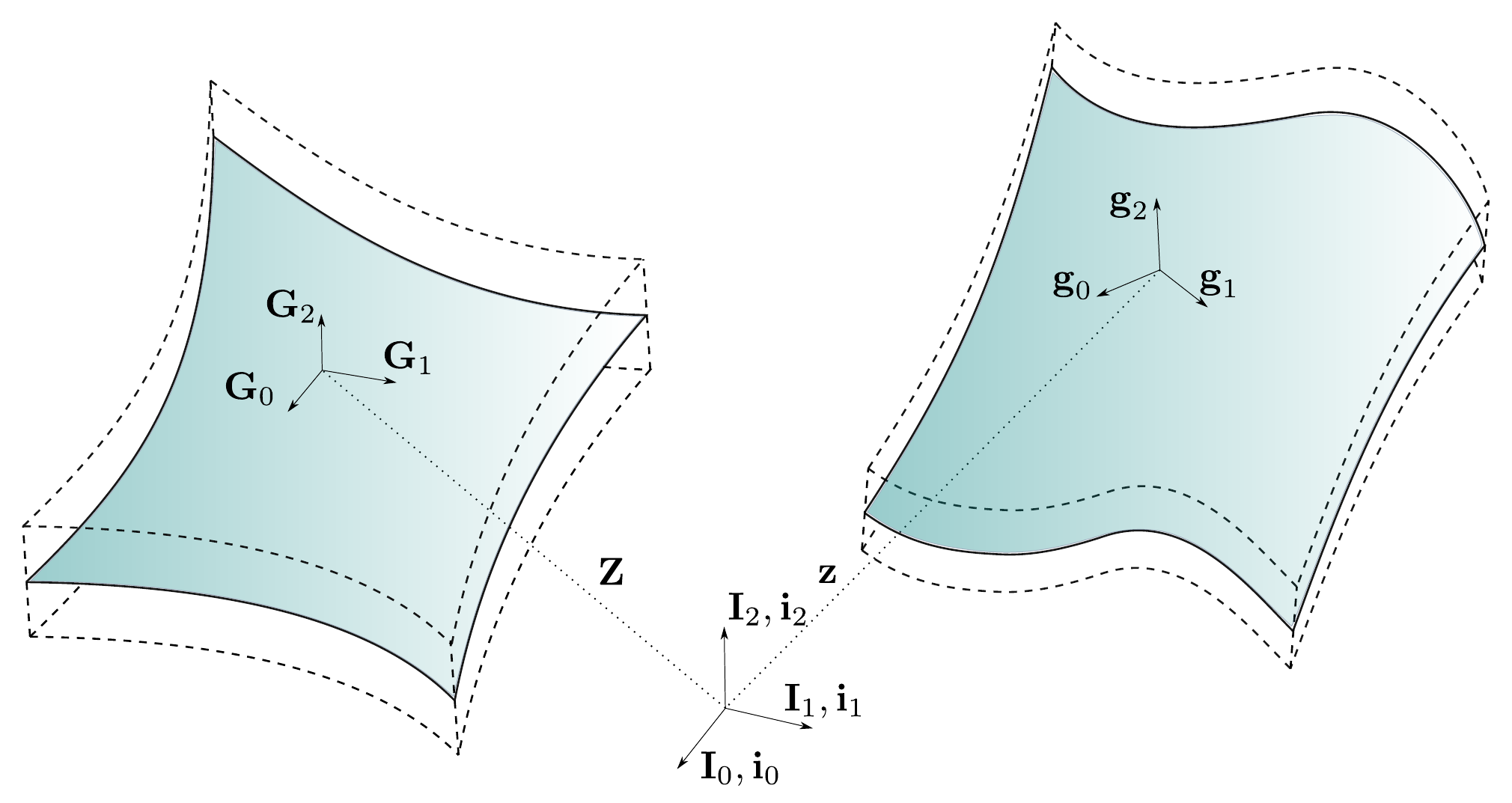}
	\caption{Shell in reference (left) and current (right) configurations.}
	\label{fig:shell_geometry}
\end{figure}

Fig. \ref{fig:shell_geometry} shows the reference (undeformed) and current (deformed)
configurations of a shell. The position $\mathbf{Z}$ of a material point in the
reference configuration can be expressed as a function of the local curvilinear
coordinates ${\vxi}=[\xi, \eta, \zeta]$ and spatially varying thickness ${}^{\mathrm{M}}a(\xi, \eta)$ as follows
\beq
\mathbf{Z}(\xi, \eta, \zeta)={}^{\mathrm{M}}\mathbf{Z}(\xi, \eta)+(\zeta-\frac{1}{2})
{}^{\mathrm{M}}a(\xi, \eta) \mathbf{D}(\xi, \eta),
\label{eqn:position_reference}
\eeq
where superscript ${}^{\mathrm{M}}(\cdot)$ indicates the field is defined only on the
shell's mid-surface. Meanwhile, ${}^{\mathrm{M}}\mathbf{Z}$ is the projection of the point
$\mathbf{Z}$ on the mid-surface and $\mathbf{D}$ represents the normal direction
vector in the reference configuration which can be defined as follows
\bseq
\begin{align}
{}^{\mathrm{M}}\mathbf{Z}(\xi, \eta) &=\frac{1}{2}\left[\mathbf{Z}_{\mathrm{u}}(\xi, \eta)+\mathbf{Z}_\ell(\xi, \eta)\right], \\ 
\mathbf{D}(\xi, \eta) &=\frac{1}{2}\left[\mathbf{Z}_{\mathrm{u}}(\xi, \eta)-\mathbf{Z}_\ell(\xi, \eta)\right], \label{eqn:projected_position_director_vector_reference}
\end{align}
\eseq	
where $\mathbf{Z}_{\mathrm{u}}$ and $\mathbf{Z}_\ell$ are the position
vectors of a point at the top and bottom surfaces in the reference
configuration, respectively. The element
kinematics is based on the linear combination of the position of the top and bottom surfaces of the element
\cite{hosseini_isogeometric_2013}. 

In the current configuration the position vector can be expressed as 
\beq
\mathbf{z}(\xi, \eta, \zeta)={}^{\mathrm{M}}\mathbf{z}(\xi, \eta)
+(\zeta-\frac{1}{2})
{}^{\mathrm{M}}a(\xi, \eta) \mathbf{d}(\xi, \eta),
\label{eqn:position_current}
\eeq
where ${}^{\mathrm{M}}\mathbf{z}$ is the projection of the point
on the mid-surface and $\mathbf{d}$ denotes the normal direction
vector in the current configuration. These are calculated as follows
\bseq
\begin{align}
{}^{\mathrm{M}}\mathbf{z}(\xi, \eta) &=\frac{1}{2}\left[\mathbf{z}_{\mathrm{u}}(\xi, \eta)+\mathbf{z}_\ell(\xi, \eta)\right], \\ 
\mathbf{d}(\xi, \eta) &=\frac{1}{2}\left[\mathbf{z}_{\mathrm{u}}(\xi, \eta)-\mathbf{z}_\ell(\xi, \eta)\right]. \label{eqn:projected_position_director_vector_current}
\end{align}
\eseq
Consequently, the position of a point in the current configuration is expressed in 
terms of its position in the reference configuration and its displacement 
$\mathbf{u}(\xi, \eta, \zeta)$ as follows
\beq
\mathbf{z}(\xi, \eta, \zeta) = \mathbf{Z}(\xi, \eta, \zeta) + \mathbf{u}(\xi, \eta, \zeta).
\eeq
	
By assuming the thickness is uniform throughout the shell in the reference configuration, the covariant base
vectors of a material point in the reference configuration can be obtained as
follows  
\bseq
\begin{align}
{\mathbf{G}}_{\alpha}&=\frac{\partial \mathbf{Z}}{\partial
\vxi^{\alpha}}=
\frac{\partial {}^{\mathrm{M}}\mathbf{Z}}{\partial \vxi^{\alpha}}
+a(\zeta-\frac{1}{2})\frac{\partial \mathbf{D}}{\partial \vxi^{\alpha}}, \quad
\alpha=0,1, \\
{\mathbf{G}}_{2}&=a\mathbf{D},
\end{align}
\eseq
where $\vxi^{\alpha}$ are basis vectors in the element's reference coordinates.

Similarly, the covariant basis vectors in the current configuration are derived
from the partial derivatives of the position vector with respect to curvilinear
coordinates $\vxi$ as follows
\bseq
\begin{align}
{\mathbf{g}}_{\alpha}&=\frac{\partial \mathbf{z}}{\partial
\vxi^{\alpha}}={\mathbf{G}}_{\alpha} + 
\frac{\partial\mathbf u}{\partial \vxi^{\alpha}}, \\
{\mathbf{g}}_{2}&={\mathbf{G}}_{2} + \frac{\partial\mathbf u}{\partial\zeta}.
\end{align}
\eseq

Consequently, the coefficients of the metric tensors in the reference and current
configurations are defined as
\bseq
\begin{align}
G^{AB} &= {\mathbf G}_A \cdot {\mathbf G}_B, \\
g_{ab} &= {\mathbf g}_a \cdot {\mathbf g}_b.
\end{align}
\eseq
The relation between the covariant and contravariant basis vectors in the reference and
current configurations are expressed as follows
\bseq
\begin{align}
	{\mathbf G}_{A} \cdot {\mathbf G}^B &= \delta_A^B, \\
	{\mathbf g}_{a} \cdot {\mathbf g}^b &= \delta_a^b,
\end{align}
\eseq
where $\delta_a^b$ is the Kronecker delta. The deformation gradient in the
curvilinear coordinates is defined as
\beq
{\mathbf F} = \frac{\partial {\mathbf z}}{\partial {\mathbf Z}} = 
\frac{\partial
{\mathbf z}}{\partial {\vxi^i}}
\frac{\partial {\vxi^i}}{\partial {\mathbf Z}}
=
\frac{\partial
{\mathbf z}}{\partial {\vxi^i}} \otimes  {\mathbf G}^i = {\mathbf g}_i \otimes {\mathbf G}^i.
\eeq

In the spirit of the solid shell formulation, any physical equation is expressed in the
current curvilinear shell coordinate system. Thus, the coefficients of the right Cauchy-Green
deformation tensor or Green's deformation tensor are expressed as
\beq
C^A_B = g_{ab}G^{AC} F^b_C F^a_B.
\eeq
Meanwhile, the Green-Lagrange strain tensor is
calculated as
\beq
E^A_B = \frac{1}{2}(C^A_B - {{\delta}_B^A}).
\eeq
For simplicity, the St.Venant-Kirchhoff material is employed and, as a consequence, the
second Piola-Kirchhoff stress tensor which relates forces in the reference
configuration with areas in the reference configuration is defined as
\beq
S^A_B = \frac{1}{2} \lambda E^A_B + \mu E^A_B E^A_B.
\eeq
On the other hand, the first Piola-Kirchhoff given which relates forces in the
current configuration with areas in the reference configuration is given by
\beq
P^A_a = g_{ac} G^{AC} F^c_B S^B_C,
\eeq
where $P^A_a$ are the coefficients of stress tensor in current local and
reference local coordinate systems. Alternately, the first Piola-Kirchhoff
stress can be defined in global coordinate system as
\beq
P^I_i = G_A^I g^a_i P^A_a.
\eeq

At this point it is worth noting that by using the proposed prism elements to discretise
shells, the thickness of the shell is not constant during deformation (i.e. in the current configuration) and
the through-thickness stretching effect is captured.
Therefore, the kinematics of the solid shell element is more comprehensive than
the kinematics assumed in the Kirchhoff-Love theory.

\section{Geometry and displacement approximations}
\label{sec:aproximations}

\subsection{Geometry approximation in the reference configuration}

In order to maintain the practicality of the approach, a surface mesh which represents the 
mid-surface of the shell is first generated using a mesh generator. 
If the mid-surface is non-planar, higher order triangle elements could 
be used to better represent its geometry using the hierarchical shape functions described in 
Section \ref{sec:prism_element}. Having the mid-surface approximated, the geometry of a shell 
structure will be fully represented by only adding its thickness. The following describes the 
approximation of the necessary components of a shell's kinematics in the reference configuration.

First, the reference position of an arbitrary point on the shell's mid-surface in the global Cartesian 
coordinate system is calculated on a triangular mesh as follow
\beq
{{}^{\mathrm{M}}Z_r}(\xi,\eta) = \sum_v \phi^v(\xi,\eta) \underline{Z}_r^v + \sum_{e_t} \sum_\ell 
\phi^{e_t}_\ell(\xi,\eta) (\underline{Z}_r)^{e_t}_\ell+ \sum_{\ell,m} \phi^t_{\ell m}(\xi,\eta) 
(\underline{Z}_r)^t_{\ell m},
\label{eqn:mid-surface_ref_position}
\eeq
where $\underline{Z}$  and the associated superscripts represent DOFs on vertex $(v)$, edge
$(e_t)$, and triangle $(t)$. Meanwhile, the right subscript $(\cdot)_r$ indicates coefficients of 
position vector in the global coordinate system. $\phi$ represents the shape functions 
associated with the DOFs on vertices, edges and triangles. The expression in Eq. 
\eqref{eqn:mid-surface_ref_position} can be written in a vector-matrix form as follow
\beq
{}^{\mathrm{M}}\mathbf{Z}(\xi,\eta) = {}^3\boldsymbol{\phi}^v(\xi,\eta) \underline{\mathbf{Z}}^v+ 
{}^3\boldsymbol{\phi}^{e_t}(\xi,\eta) \underline{\mathbf{Z}}^{e_t}+ {}^3\boldsymbol{\phi}^{t}(\xi,\eta) 
\underline{\mathbf{Z}}^{t}= \sum_{g=v,e_t,t} {}^3\boldsymbol{\phi}^g(\xi,\eta) 
\underline{\mathbf{Z}}^{g},
\eeq
where the left superscript $^3(\cdot)$ indicates the matrix of shape functions for a 
three-component vector field.

Second, the through-thickness direction vector, of unit length, of a point on the mid-surface is 
approximated as follow
\beq
{}^{\mathrm{M}}\mathbf{D}(\xi,\eta)= \sum_{g=v,e_t,t} {}^3\boldsymbol{\phi}^g(\xi,\eta) \underline{\mathbf{D}}^{g}.
\eeq

Third, the shell's thickness at a point on the mid-surface is calculated as follow
\beq
{}^{\mathrm{M}}a(\xi,\eta)= \sum_{g=v,e_t,t} {}^1\boldsymbol{\phi}^g(\xi,\eta) \underline{\mathbf{a}}^{g},
\eeq
where the left superscript $^1(\cdot)$ denotes the matrix of shape functions for a scalar field.

Finally, having all the necessary fields described and approximated, the reference position in the global 
Cartesian coordinate system of a point in a prism element is given by
\beq
\mathbf{Z}(\xi,\eta,\zeta) = {}^{\mathrm{M}}\mathbf{Z}(\xi,\eta)+ 
{}^{\mathrm{M}}a(\xi,\eta)(\zeta-\frac{1}{2}){}^{\mathrm{M}}\mathbf{D}(\xi,\eta).
\label{eqn:shell_ref_position}
\eeq
Similarly, it can be expressed in a compact form as follow 
\beq
\mathbf{Z}(\xi,\eta,\zeta) = \sum_{s=L,U}\sum_{g=v,e_t,t} {}^3\boldsymbol{\phi}^g_s(\xi,\eta,\zeta) 
\underline{\mathbf{Z}}^{g}_{s}.
\eeq

Next, the Jacobian is derived to transform a point from a prism element's local coordinate system 
to the global Cartesian coordinate system. In order 
to keep the derivation 
straightforward, the thickness of the shell is assumed constant, i.e. 
${}^{\mathrm{M}}a(\xi,\eta)=\textrm{const}=a$. Using Eq. \eqref{eqn:shell_ref_position}, the Jacobian is then 
defined as follow
\beq
\mathbf{J}(\xi,\eta,\zeta) = \frac{\partial \mathbf{Z}}{\partial \boldsymbol{\xi}} = 
{}^{\mathrm{M}}\mathbf{J}(\xi,\eta) + {}^D\mathbf{J}(\xi,\eta,\zeta),
\eeq
where the first term contributed by the mid-surface components is
\beq
{}^{\mathrm{M}}\mathbf{J}(\xi,\eta,\zeta) = \left[ \begin{array}{ccc} \frac{\partial {}^{\mathrm{M}}_0Z}{\partial \xi} & 
\frac{\partial {}^{\mathrm{M}}_0Z}{\partial \eta} & 0 \\ \frac{\partial {}^{\mathrm{M}}_1Z}{\partial \xi} & \frac{\partial 
{}^{\mathrm{M}}_1Z}{\partial \eta} & 0 \\ 0 & 0 & 0 \end{array} \right],
\eeq
and the second term contributed by the through-thickness component is
\beq
{}^D\mathbf{J}(\xi,\eta,\zeta) = a \left[ \begin{array}{ccc} 0 & 0 & {}^{\mathrm{M}}_0D \\ 0 & 0 & {}^{\mathrm{M}}_1D \\ 0 & 
0 & {}^{\mathrm{M}}_2D \end{array} \right] + a(\zeta-\frac{1}{2}) \left[ \begin{array}{ccc} \frac{\partial 
{}^{\mathrm{M}}_0D}{\partial \xi} & \frac{\partial {}^{\mathrm{M}}_0D}{\partial \eta} & 0 \\ \frac{\partial {}^{\mathrm{M}}_1D}{\partial \xi} & 
\frac{\partial {}^{\mathrm{M}}_1D}{\partial \eta} & 0 \\ \frac{\partial {}^{\mathrm{M}}_2D}{\partial \xi} & \frac{\partial 
{}^{\mathrm{M}}_2D}{\partial \eta} & 0 \end{array} \right].
\eeq

\subsection{Curvilinear systems in reference and current configuration}

\begin{figure}[!htbp]
	\centering
	\includegraphics[scale=0.3]{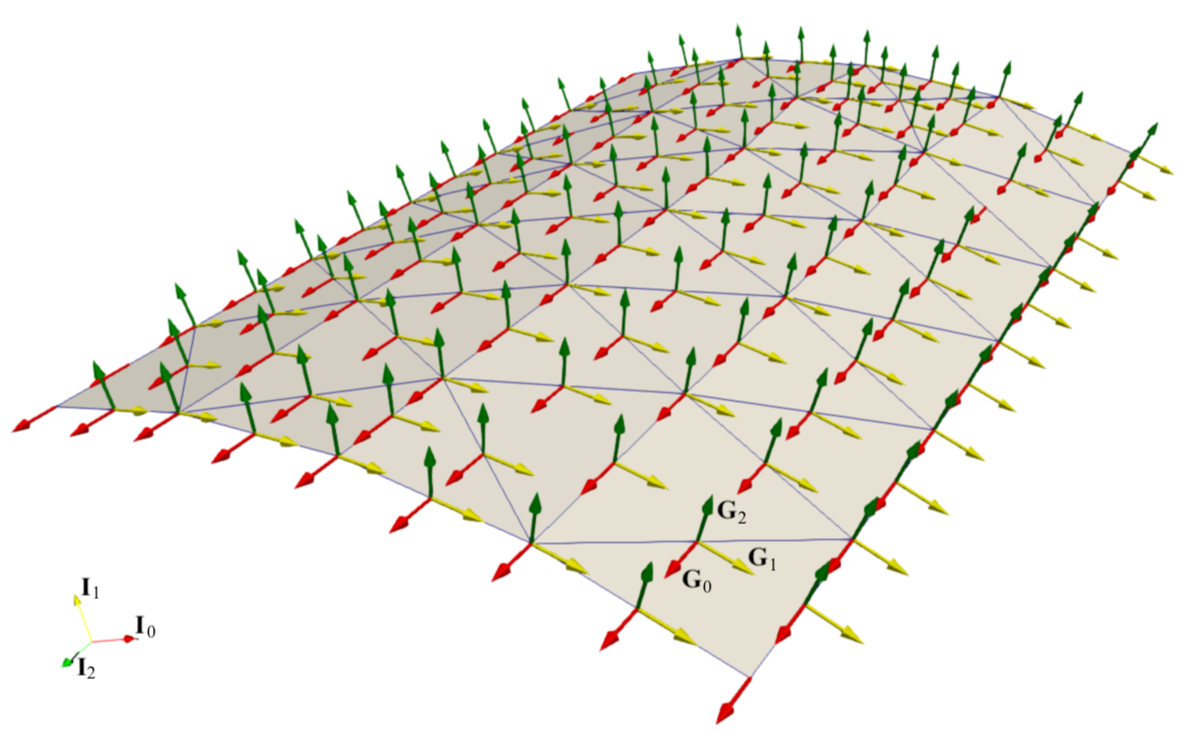}
	\caption{Approximated local curvilinear basis vector $\mathbf{G}_0$ (red), $\mathbf{G}_1$ (yellow), and $\mathbf{G}_2$ (green) on the mesh of a curved shell.}
	\label{fig:approximated_curvilinear_base_vectors}
\end{figure}

A local curvilinear basis in the reference configuration, with basis vectors ${\bf G}_\alpha$, is established. The in-plane components of ${\bf G}_\alpha$ are
approximated using hierarchical approximation functions as follows
\bseq
\begin{align}
	{}^{\mathrm{M}}\mathbf{G}_0(\xi,\eta) &= G^0_I \mathbf{I}^I = \left\{ \sum_{g=v,e_t,t} {}^3\boldsymbol{\phi}^g(\xi,\eta) \underline{\mathbf{G_0}}^{g} \right\} \mathbf{I}^I, \\
	{}^{\mathrm{M}}\mathbf{G}_1(\xi,\eta) &= G^1_I \mathbf{I}^I = \left\{ \sum_{g=v,e_t,t} {}^3\boldsymbol{\phi}^g(\xi,\eta) \underline{\mathbf{G_1}}^{g} \right\} \mathbf{I}^I,
\end{align}
\eseq
where ${\bf I}$ is the global Cartesian basis vector, as shown in Fig. \ref{fig:shell_geometry}. Meanwhile, the local basis vector in the through-thickness direction is given by
\beq
{}^{\mathrm{M}}\mathbf{G}_2(\xi,\eta) = G^2_I \mathbf{I}^I = \left\{
\textrm{Spin}[{}^{\mathrm{M}}\mathbf{G}_0(\xi,\eta)]{}^{\mathrm{M}}\mathbf{G}_1(\xi,\eta)
\right\} \mathbf{I}^I,
\label{eqn:base_vector_thickness}
\eeq
where $\textrm{Spin}[\cdot]$ is the spin operator acting as a vector product.
In  particular, for an axis ${a_i}$, the spin tensor is defined as
\beq
\operatorname{Spin}[{\bf a}]=\left[\begin{array}{ccc}0 & -a_{3} & a_{2} \\ a_{3} 
& 0 & -a_{1} \\ -a_{2} & a_{1} & 0\end{array}\right],
\eeq
and the cross product of ${\bf a}$ and an arbitrary vector ${\bf b}$ is expressed as
\beq
\mathbf{a} \times \mathbf{b}=\operatorname{Spin}[\mathbf{a}]
\mathbf{b}=-\operatorname{Spin}[\mathbf{b}] \mathbf{a}=-\mathbf{b} \times \mathbf{a}.
\eeq
The local curvilinear basis vectors
$\mathbf{G}_0, \mathbf{G}_1,$ and $\mathbf{G}_2$ are visualised in Fig.
\ref{fig:approximated_curvilinear_base_vectors}.

The local curvilinear coordinate basis in the current configuration, which is convected by shell mid-surface motion, is given as 
\beq
\mathbf{g}_a = g_a^i \mathbf{i}_i = {}^{\mathrm{M}}F^i_I\,\, {}^{\mathrm{M}}G^{AI} \delta^a_A \mathbf{i}_i,
\eeq
where ${}^{\mathrm{M}}F^i_I(\xi,\eta,\zeta)$ are local components of the deformation gradient.
It should be noted that, for a discretised system, the deformation gradient is 
additively decomposed into a part resulting from through-thickness DOFs and a part from the upper and lower triangles' DOFs. 
Basis vector $\mathbf{i}_i$ is a Cartesian basis vector for the current  
configuration, which could be set independently from the Cartesian basis 
$\mathbf{I}_i$. However, for the rest of paper, we restrict ourself to the case 
where $\mathbf{i}_i \equiv \mathbf{I}_i$. 
This can be considered a generalisation of the co-rotational
formulation for shells.

\subsection{Displacement approximations}
\begin{figure}[!htbp]
	\centering
	\includegraphics[scale=.7]{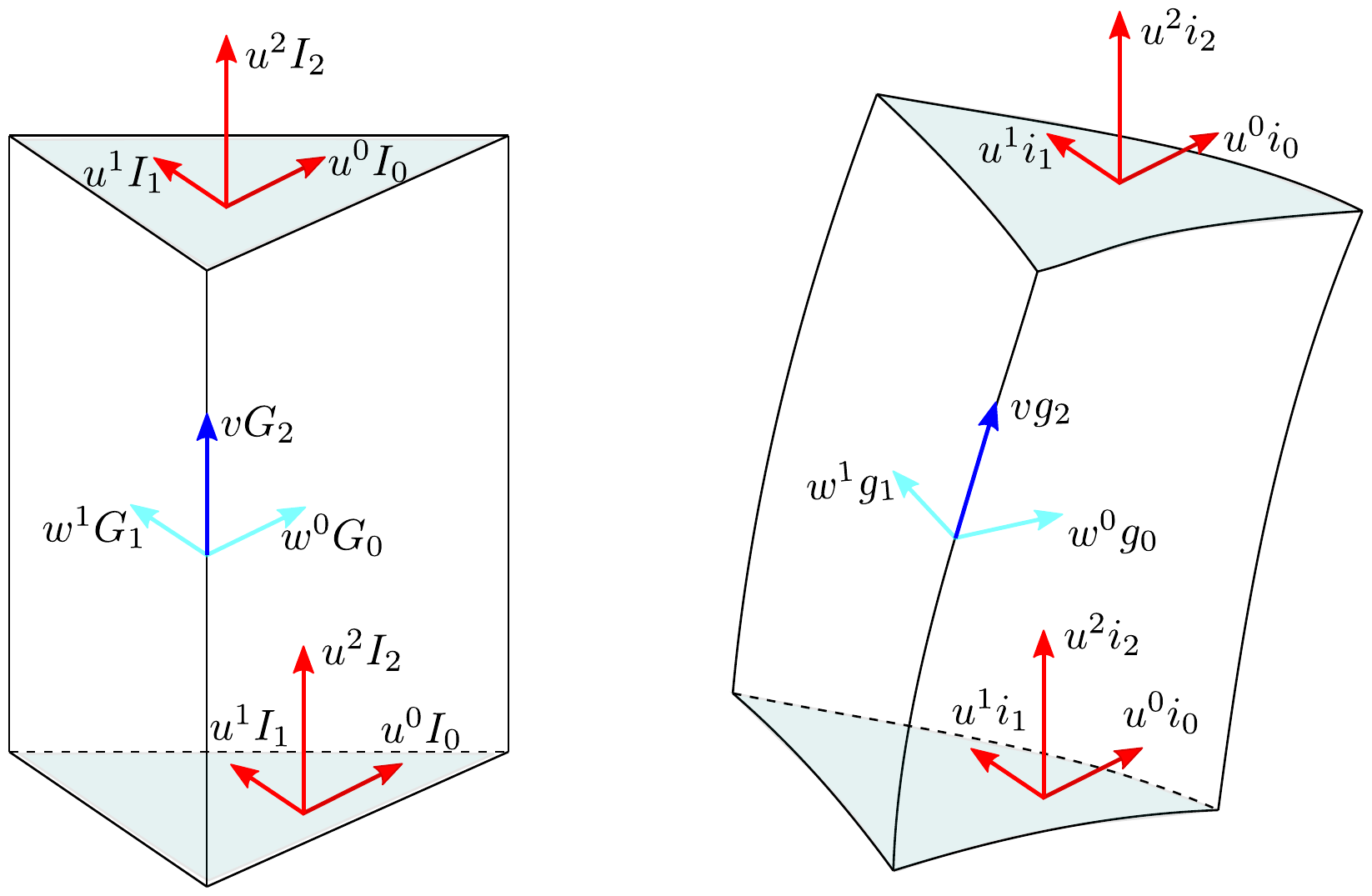}
	\caption{Displacement components ${\bf u}$, ${\bf v}$, and ${\bf w}$
	and their associated bases (global Cartesian bases ${\bf I, i}$ and local
	convected curvilinear bases ${\bf G, g}$) in the reference (left) and current
	(right) configurations.}
	\label{fig:displacement_components}
\end{figure}

The position of a material point in the current configuration is given by 
\beq
z^J\mathbf{I}_J = Z^J\mathbf{I}_J + u^J \mathbf{I}_J + v g_2^J \mathbf{I}_J
 + w^\alpha g_\alpha^J \mathbf{I}_J,
\label{eqn:position_current_config}
\eeq
where the displacement coefficients $u^J$ are DOFs associated with the vertices, edges and faces of the upper and lower triangles, 
and expressed in the global Cartesian coordinate system.
The displacement coefficients $v$ and $w^\alpha$ are DOFs associated with the edges and faces of the quadrilaterals through the shell thickness, and
expressed in terms of the current curvilinear (convected) coordinate system. 
These displacement DOFs are
illustrated in Fig. \ref{fig:displacement_components}.

The displacement vector ${\bf u}$
on the top and bottom triangles is approximated as follow
\beq
{\bf u}(\boldsymbol \xi) = \sum_{s=L,U} \sum_{g=v,e_t,t} {}^3\boldsymbol{\phi}^g_s(\boldsymbol \xi) \underline{\mathbf{u}}^{g}_s = \left[{}^3\boldsymbol{\overline{\phi}}(\boldsymbol \xi)\right] \underline{\mathbf{u}},
\eeq
with the displacement DOFs in vector
$\underline{\mathbf{u}}$ expressed in the global Cartesian coordinate system. This allows
elements to be assembled with other tetrahedral elements without the need for
any linking or transfer elements. 
Meanwhile, the through-thickness displacements
$\mathbf{v}=v g_2^J \mathbf{I}_J$ and $\mathbf{w}=w^\alpha g_\alpha^J \mathbf{I}_J$
are approximated in the local curvilinear system which follows the global
mid-surface deformation. The coefficients $v$ and $w^\alpha$ are given by 
\bseq
\begin{align}
v(\boldsymbol \xi) &= \sum_{g=e_q,q,p} {}^1\boldsymbol{\phi}^g(\boldsymbol \xi) \underline{\mathbf{v}}^{g} = \left[{}^1\boldsymbol{\overline{\phi}}(\boldsymbol \xi)\right]\underline{\mathbf{v}}, \\
\mathbf{w}(\boldsymbol \xi) &= \sum_{g=e_q,q,p} {}^2\boldsymbol{\phi}^g(\boldsymbol \xi) \underline{\mathbf{w}}^{g} = \left[{}^2\boldsymbol{\overline{\phi}}(\boldsymbol \xi)\right]\underline{\mathbf{w}}.
\end{align}
\eseq


Having ${\bf u}$, ${\bf v}$ and ${\bf w}$ derived, the gradient of displacements can be calculated as follows
\bseq
\begin{align}
	\nabla_{Z^J} u^i \left( \mathbf{I}_i \otimes \mathbf{I}^J \right) &= \left\{ \sum_{s=U,L}\sum_{g=v,e_t,t} \frac{\partial {}^3 \{\boldsymbol \phi^g_s\}^i \mathbf{I}_i }{\partial \xi^k} J_J^k \right\} \underline{\mathbf{u}}^g_s \otimes \mathbf{I}^J, \\
	\nabla_{X^A} u^a &= g_j^a \nabla_{Z^J} u^j G^J_A,
\end{align}
\eseq
and the gradients of higher-order through thickness displacements
\bseq
\begin{align}
\nabla_{X^A} v^a \left( \mathbf{g}_a \otimes \mathbf{G}^A \right) = \left\{ \sum_{g=e_q,q,p} \frac{\partial {}^1 \{\boldsymbol \phi^g\}^a \mathbf{g}_a }{\partial \xi^k} J_J^k G_A^J \right\} \underline{v}^{g} \otimes \mathbf{G}^A, \\
\nabla_{X^A} w^a \left( \mathbf{g}_a \otimes \mathbf{G}^A \right) = \left\{ \sum_{g=e_q,q,p} \frac{\partial {}^2 \{\boldsymbol \phi^g\}^a \mathbf{g}_a }{\partial \xi^k} J_J^k G_A^J \right\} \underline{w}^{g} \otimes \mathbf{G}^A.
\end{align}
\eseq
Consequently, the deformation gradient is calculated as
\bseq
\begin{align}
	{}^{\mathrm{M}}F^i_I &= \delta^i_I + \nabla_{Z^I} u^i = \delta^i_I + \left[\nabla_{Z^I}{}^3\boldsymbol{{\phi}}^i(\boldsymbol \xi)\right] \underline{\mathbf{u}}, \\
	F^a_A &= g^a_iG^I_A\delta^i_I  + \nabla_{X^A} v^a + \nabla_{X^A} w^a= g^a_i G^I_A \delta^i_I + \left[\nabla_{X^A}{}^1\boldsymbol{{\phi}}(\boldsymbol \xi)\right] \underline{\mathbf{v}} + \left[\nabla_{X^A}{}^2\boldsymbol{{\phi}}(\boldsymbol \xi)\right] \underline{\mathbf{w}}.
\end{align}
\eseq
It will be shown later in Section~\ref{sec:example} that the approximation order for the in-plane displacements through the thickness of the shell should be at
least second-order to avoid locking.
It should be noted that when only DOFs are set on triangle entities, the
approximation of displacement through the thickness is linear. The DOFs
associated with the quadrilaterals' edges and faces, and prism volume are higher orders and
defined in curvilinear coordinate system which is convected by the deformation
gradient ${}^{\mathrm{M}}{\bf F}$.

\subsection{Element formulation}

With all the physical quantities derived and approximated, the vector of 
internal forces for the top/bottom surfaces and though-thickness direction are
respectively expressed as follows
\bseq
\begin{align}
	\mathbf{f}^\textrm{int}_\textrm{u}(\underline{\mathbf{u}},\underline{\mathbf{v}},\underline{\mathbf{w}}) &= \int_\Omega \left[\nabla_{Z^I}{}^3\boldsymbol{\overline{\phi}} \right] P^I_i(\underline{\mathbf{u}},\underline{\mathbf{v}},\underline{\mathbf{w}}) \textrm{d}\Omega, \\
	\mathbf{f}^\textrm{int}_\textrm{v}(\underline{\mathbf{u}},\underline{\mathbf{v}},\underline{\mathbf{w}}) &= \int_\Omega \left[\nabla_{X^A}{}^1\boldsymbol{\overline{\phi}} \right] P^A_a(\underline{\mathbf{u}},\underline{\mathbf{v}},\underline{\mathbf{w}}) \textrm{d}\Omega, \\
	\mathbf{f}^\textrm{int}_\textrm{w}(\underline{\mathbf{u}},\underline{\mathbf{v}},\underline{\mathbf{w}}) &= \int_\Omega \left[\nabla_{X^A}{}^2\boldsymbol{\overline{\phi}} \right] P^A_a(\underline{\mathbf{u}},\underline{\mathbf{v}},\underline{\mathbf{w}}) \textrm{d}\Omega, 
\end{align}
\eseq
where $\Omega$ is the domain occupied by the shell.

The nonlinear problem is solved using the Newton-Raphson method in which a Taylor series expansion of the residual vector is truncated after the linear term:
\beq
\left[ \begin{array}{ccc} 
  \frac{\partial \mathbf{f}^\textrm{int}_\textrm{u}(\underline{\mathbf{u}},\underline{\mathbf{v}},\underline{\mathbf{w}})}{\partial\underline{\mathbf{u}}} & 
  \frac{\partial \mathbf{f}^\textrm{int}_\textrm{u}(\underline{\mathbf{u}},\underline{\mathbf{v}},\underline{\mathbf{w}})}{\partial\underline{\mathbf{v}}} & 
  \frac{\partial \mathbf{f}^\textrm{int}_\textrm{u}(\underline{\mathbf{u}},\underline{\mathbf{v}},\underline{\mathbf{w}})}{\partial\underline{\mathbf{w}}} \\
   \frac{\partial \mathbf{f}^\textrm{int}_\textrm{v}(\underline{\mathbf{u}},\underline{\mathbf{v}},\underline{\mathbf{w}})}{\partial\underline{\mathbf{u}}} &
   \frac{\partial \mathbf{f}^\textrm{int}_\textrm{v}(\underline{\mathbf{u}},\underline{\mathbf{v}},\underline{\mathbf{w}})}{\partial\underline{\mathbf{v}}} &
   \frac{\partial \mathbf{f}^\textrm{int}_\textrm{v}(\underline{\mathbf{u}},\underline{\mathbf{v}},\underline{\mathbf{w}})}{\partial\underline{\mathbf{w}}} \\
   \frac{\partial \mathbf{f}^\textrm{int}_\textrm{w}(\underline{\mathbf{u}},\underline{\mathbf{v}},\underline{\mathbf{w}})}{\partial\underline{\mathbf{u}}} &
   \frac{\partial \mathbf{f}^\textrm{int}_\textrm{w}(\underline{\mathbf{u}},\underline{\mathbf{v}},\underline{\mathbf{w}})}{\partial\underline{\mathbf{v}}} &
   \frac{\partial \mathbf{f}^\textrm{int}_\textrm{w}(\underline{\mathbf{u}},\underline{\mathbf{v}},\underline{\mathbf{w}})}{\partial\underline{\mathbf{w}}} \end{array} \right] 
  \left\{ \begin{array}{c} \delta \underline{\mathbf{u}} \\ \delta \underline{\mathbf{v}} \\ \delta \underline{\mathbf{w}} \end{array} \right\} 
  = \left[ \begin{array}{c} \mathbf{f}^\textrm{ext}_\textrm{u} -
  \mathbf{f}^\textrm{int}_\textrm{u}(\underline{\mathbf{u}},\underline{\mathbf{v}},\underline{\mathbf{w}})\\
  -
  \mathbf{f}^\textrm{int}_\textrm{v}(\underline{\mathbf{u}},\underline{\mathbf{v}},\underline{\mathbf{w}})
  \\ -
  \mathbf{f}^\textrm{int}_\textrm{w}(\underline{\mathbf{u}},\underline{\mathbf{v}},\underline{\mathbf{w}})
  \end{array} \right], \label{eqn:system_equations}
\eeq
where $\delta \underline{\mathbf{u}}$, $\delta \underline{\mathbf{v}}$ and $\delta \underline{\mathbf{w}}$ are
the iterative changes in displacements in the global and current local coordinate systems, and the left hand matrix represent the tangent matrix.
Additionally, since the solution is strongly nonlinear and potentially includes
local instabilities, the classical quadratic arc-length method with line
search is employed.

As can be seen in Eq. \eqref{eqn:system_equations}, the external forces are applied only to the 
DOFs associated with the upper and lower triangles, which are expressed in the global Cartesian coordinate system. Meanwhile,
all DOFs associated with the convected curvilinear system are internal and there is no
external forces act on them, except body forces.

It is worth noting that in this proposed formulation, the current curvilinear
system follows the shell deformation described by DOFs on the upper and lower triangle. 
This is similar to the co-rotational formulation presented by Crisfield
\cite{crisfield_consistent_1990} in which the coordinate system is rotated while
the shell is deforming and the tangent stiffness in non-symmetric. However, for
the present problem, numerical experiments show that the tangent stiffness
matrix becomes symmetric as the iterative procedure reaches equilibrium state.
Therefore, the matrix can be symmetrised without deteriorating the convergence
order. 

\section{Numerical examples}
\label{sec:example}

\begin{table}[!htbp]
	\centering
	\caption{Geometry, material, and loading data}
	\begin{tabular}{llllllll}
		\hline
		Parameter & \multicolumn{2}{l}{Linear analysis$^\dagger$} & & \multicolumn{4}{l}{Nonlinear analysis$^\ddagger$} \\
		\cline{2-3}  \cline{5-8}   & PC-L       & SLR &  &  SAP &
		HS & POC & PC \\ 
		\hline
		\emph{Geometry} &   &       &       &      &     &       &  \\
		Length, $L$ & 600  & 50    &    & N/A   & N/A   & 10.35 & 200 \\
		Radii, $R$/$R_i$/$R_o$ & 300  & 25     &   & 6/10  & 10    & 4.953 & 100 \\
		Thickness, $a$ & 3  & 0.25    &  & 0.03  & 0.04  & 0.094 & 1 \\
		Angle, $\theta$ & N/A  & $40^0$   &  & N/A   & $18^0$ & N/A   & N/A \\
		&       &   &       &       &        &   &  \\
		\emph{Material} &   &       &        &   &       &       &  \\
		Young's modulus,  $E$ & 3  & $4.32\times 10^8$   &  & $21\times 10^6$ & $6.825
		\times 10^7$ & $10.5 \times 10^6$ & $30 \times 10^3$ \\
		Poisson's ratio, $\nu$ & 0.3  & 0.0  &     & 0.0   & 0.3   & 0.3125 & 0.3 \\
		&       &       &    &      &    &       &  \\
		\emph{Loading} &    &      &       &     &      &       &  \\
		Point force, $P$ &  1 & N/A  &     & N/A   & 400   & 40,000 & 12,000 \\
		Distributed force, $q$ & N/A  & 90   &     & 0.8   & N/A   & N/A   & N/A \\
		\hline
		\multicolumn{8}{l}{$^\dagger$PC-L - Pinched cylinder for linear analysis;
		SLR - Scordelis-Lo Roof}\\
		\multicolumn{8}{l}{$^\ddagger$ SAP - Slit annular
			plate;  HS - Hemispherial shell; POC - Pullout cylinder; PC - Pinched cylinder}
	\end{tabular}
	\label{tab:input_data}
\end{table}

\subsection{Linear analysis}

\subsubsection{Pinched cylinder}
\begin{figure}
	\hspace{-1cm}
	~\input{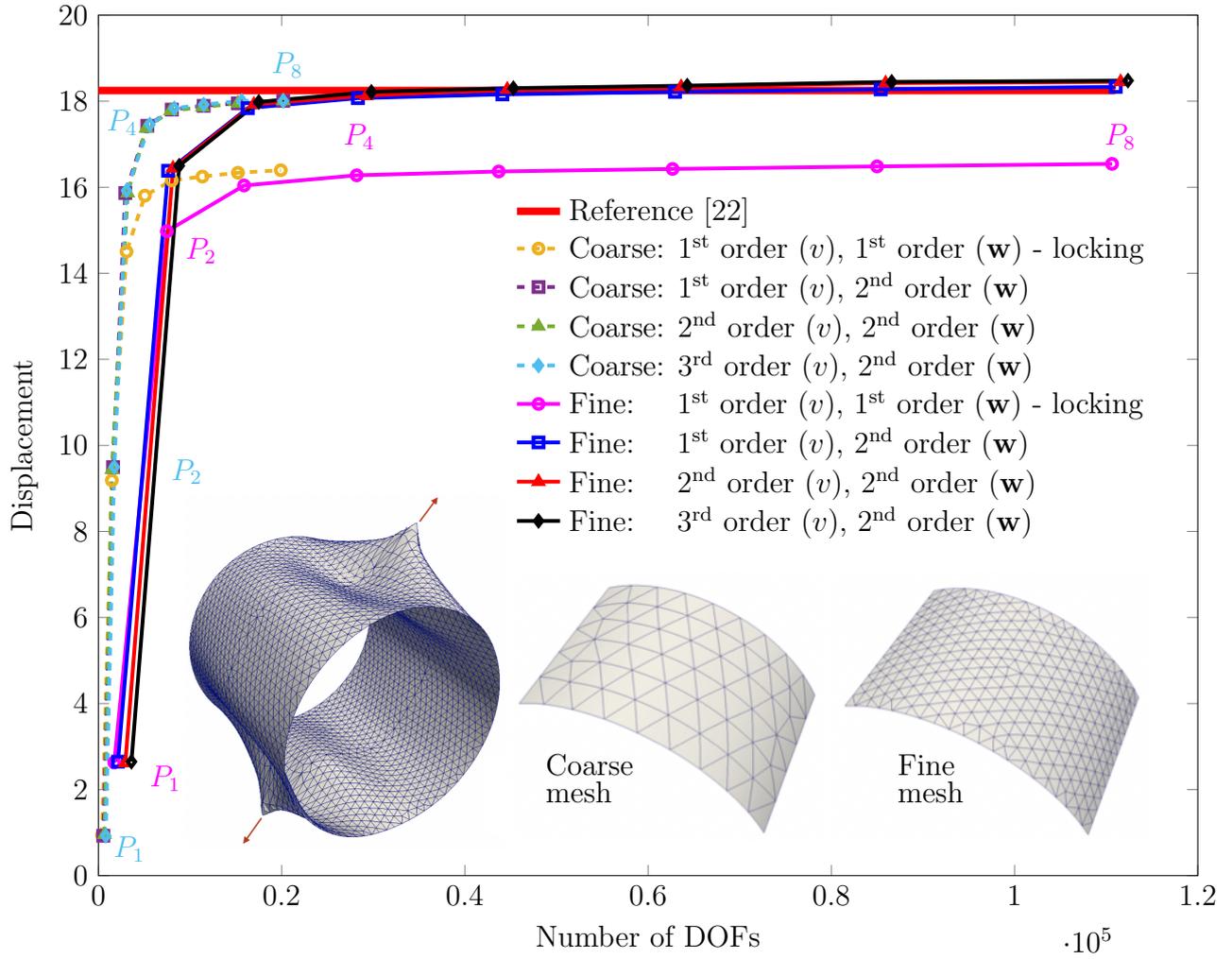}
	\caption{Convergence of displacement with different meshes and polynomial orders. $P_n$ indicates result using $n^{\rm th}$ order polynomial for displacement component ${\bf u}$ (applies to all lines).} 
	\label{fig:linear_pinched_cylinder}
\end{figure}

A linear analysis of a cylinder mounted on a diaphragm, subject to pinching forces,
is presented to demonstrate the accuracy of the proposed
approach as well as its capability to avoid locking effects. The input
data can be found in Table \ref{tab:input_data}. Only one eighth
of the geometry is discretised due to symmetry. As can be seen
in Fig. \ref{fig:linear_pinched_cylinder}, both meshes (coarse and fine) yield good convergence results compared to those reported in the
literature \cite{belytschko_stress_1985}. 
For both meshes, there are clear
signs of locking when only linear polynomials are employed for the through-thickness approximation of the in-plane displacement, $\mathbf{w}$. 
However, when this is increased to second order, the locking phenomenon is alleviated for 
both meshes. 

Fig.~\ref{fig:linear_pinched_cylinder} shows convergence plots for both the coarse and fine meshes as the order of 
approximation for the trianglular faces of the prism elements are increased from first to eighth order.
It can also be observed that only increasing the through-thickness order of approximation for the 
 displacement in the through-thickness direction $v$ has minimal effects on the results. 
Therefore, second-order 
approximations will also be applied for $v$ for
the rest of the examples presented in this paper, unless otherwise stated. It should be noted
that for visualisation purposes and because of the linear nature of the problem,
the deformed shape in Fig. \ref{fig:linear_pinched_cylinder} is presented as if
the cylinder is subject to pullout forces rather than pinched, as shown in the
literature. Since the small displacement assumption is made in this example for 
comparison purpose, the direction of the load, inward or outward, does not
change the character of the results.

\subsubsection{Scordelis-Lo Roof}
\label{sec:scordelis_lo_roof}

In this example, analysis of the well-known Scordelis-Lo Roof is considered in 
order to demonstrate both the performance of the prism element and how
the properties of the heterogeneous and hierarchical shape functions can be 
utilised to build an efficient multi-grid solver.

\begin{figure}[!htbp]
	\centering
	\def\svgwidth{8cm}
	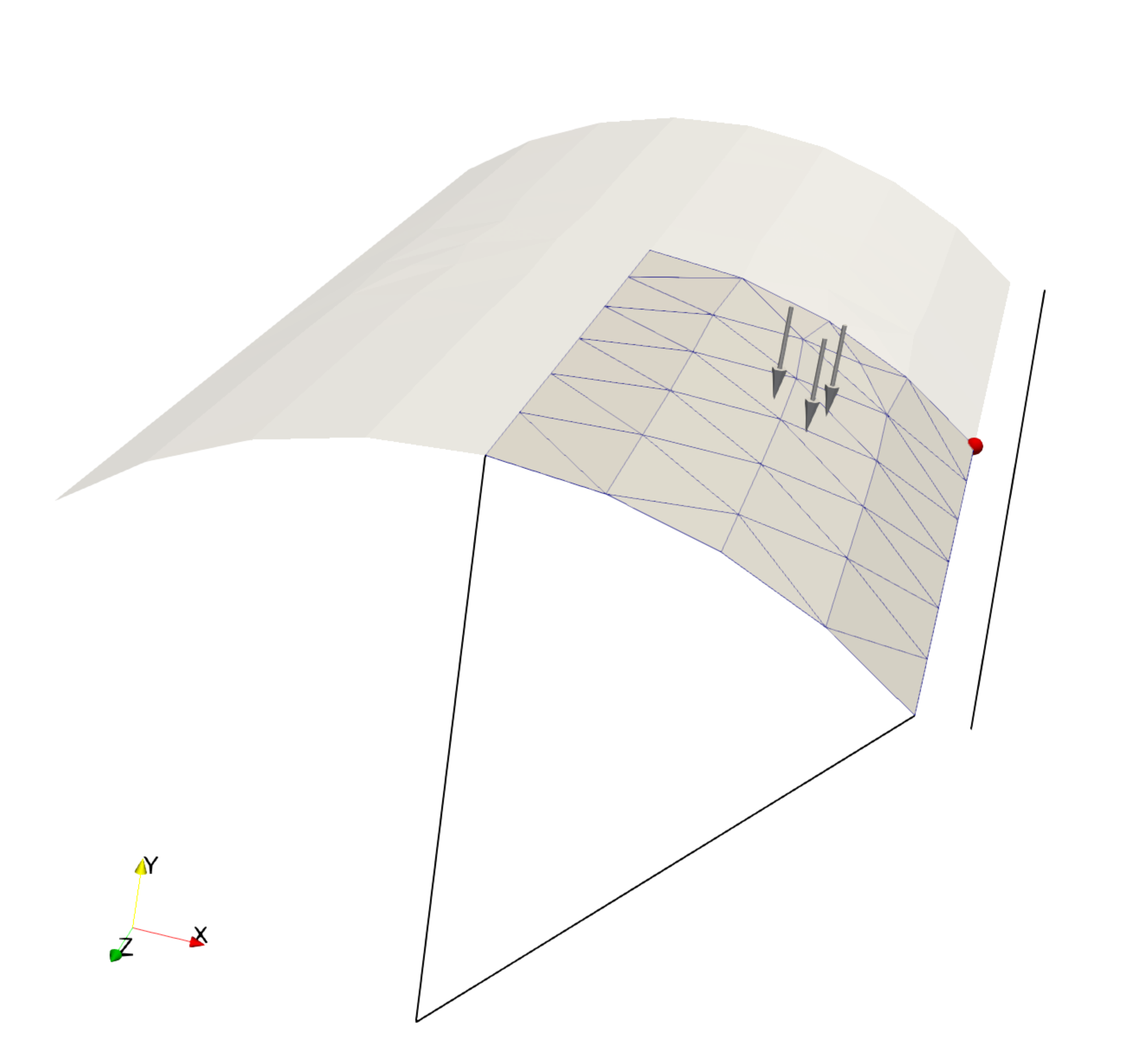	
	\caption{Scordelis-Lo roof geometry, boundary conditions and loading.}
	\label{fig:scordelis_lo_roof_geometry}
\end{figure}

Fig. \ref{fig:scordelis_lo_roof_geometry} shows the configuration of the roof
which is clamped at two ends and free along the other edges. 
The input geometry and material parameters are provided in Table \ref{tab:input_data}.
The roof is subjected to a uniformly distributed self-weight load of
$q=90$ per unit area. Due to symmetry, only one
quadrant of the full roof is considered for analysis as shown in Fig. \ref{fig:scordelis_lo_roof_geometry}.

\begin{figure}[!htbp]
	\centering
	\includegraphics[scale=0.13]{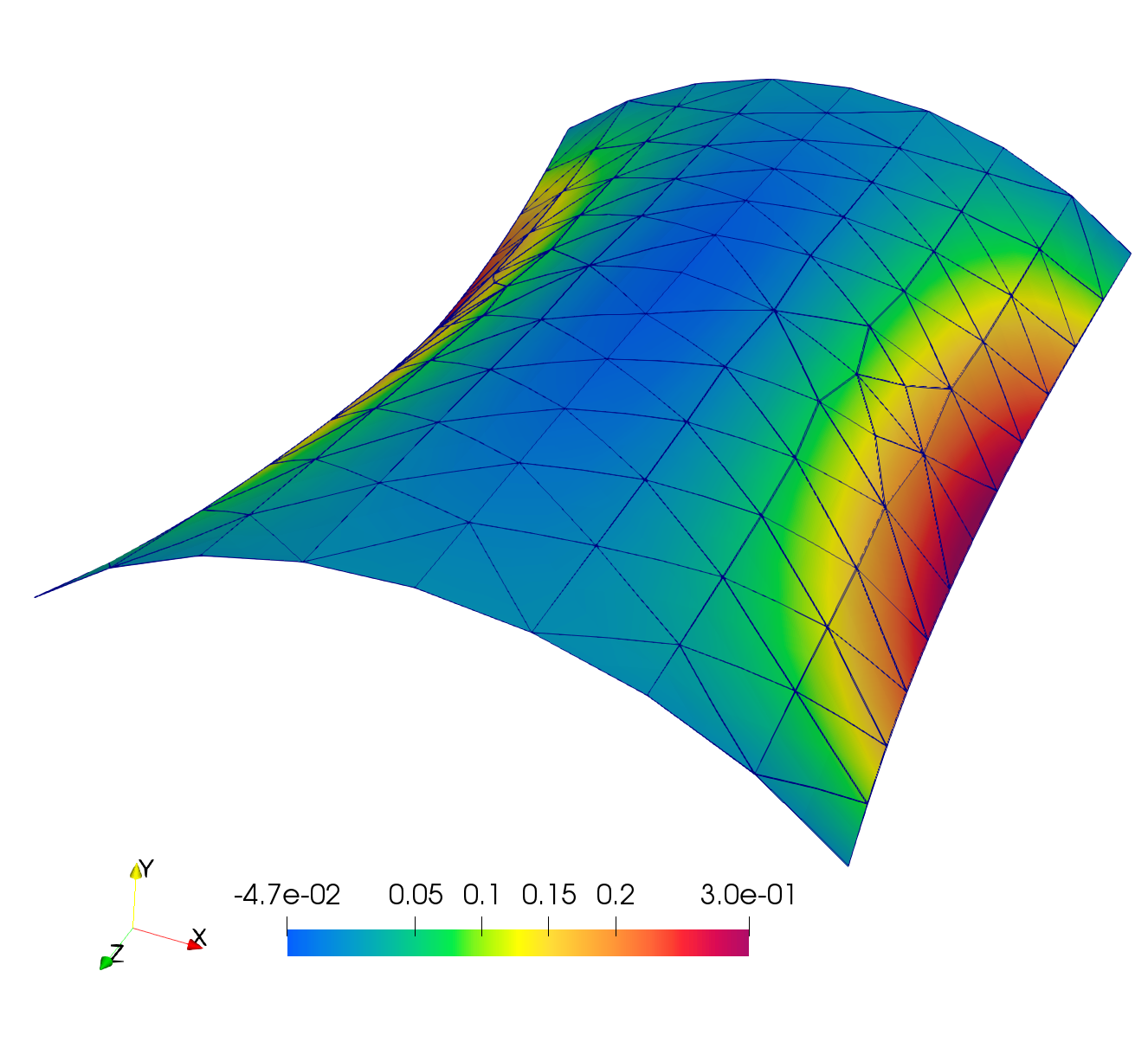}
	\caption{Deformation of Scordelis-Lo roof with fifth-order approximations for the upper and lower faces. Contours of vertical displacement.}
	\label{fig:scordelis_lo_roof_displacement}
\end{figure}

As validation, a linear analysis of the roof was undertaken and the vertical
displacement at point A (middle of the free edge as shown in Fig.
\ref{fig:scordelis_lo_roof_geometry}) is observed and compared 
to a reference value provided by Belytschko et al. \cite{belytschko_stress_1985}. 
By employing elements with anisotropic
shape functions, in which fifth order is used for the approximation functions on the upper and lower faces and second-order through the thickness, 
the vertical displacement at point A is calculated as 0.3176. 
This is slightly softer but comparable to the reference results of 
0.3024~\cite{belytschko_stress_1985} in which a 9-node
element with uniform $2\times2$ quadrature and 
$\gamma$-stabilisation was used.
The deformed shape and vertical displacement of the full structure of
Scordelis-Lo Roof is given in Fig. \ref{fig:scordelis_lo_roof_displacement}.

\begin{figure}[!htb]
	\centering
	~
%
%
\begin{tikzpicture}

\begin{axis}[%
width=4in,
height=2.8in,
at={(1.011in,0.642in)},
scale only axis,
xmode=log,
xmin=1000,
xmax=1000000,
xminorticks=true,
xlabel style={font=\color{white!15!black}},
xlabel={Number of DOFs},
ymode=log,
ymin=1e-08,
ymax=0.001,
yminorticks=true,
ylabel style={font=\color{white!15!black}},
ylabel={Error (energy norm)},
axis background/.style={fill=white},
legend style={legend cell align=left, align=left, fill=none, draw=none}
]
\addplot [color=blue, line width=1.5pt, mark=o, mark options={solid, blue}]
  table[row sep=crcr]{%
6856.919206	0.000734635\\
8568.468921	0.000450835\\
10424.72281	0.000298845\\
12022.64435	0.000174209\\
14497.4067	7.27122e-05\\
17020.32082	3.03489e-05\\
19982.28557	1.4404e-05\\
22437.15261	9.30572e-06\\
};
\addlegendentry{Coarse mesh}

\addplot [color=green, line width=1.5pt, mark=square, mark options={solid, green}]
  table[row sep=crcr]{%
15708.29441	0.000512651\\
19982.28557	0.000208543\\
24311.20311	7.65471e-05\\
28541.96511	2.47091e-05\\
33808.99476	7.01426e-06\\
39692.61186	2.44564e-06\\
44568.93516	1.89141e-06\\
49600.25063	1.84342e-06\\
};
\addlegendentry{Middle mesh}

\addplot [color=red, line width=1.5pt, mark=triangle, mark options={solid, red}]
  table[row sep=crcr]{%
54709.72303	0.000112548\\
68977.85379	2.28757e-05\\
83921.06135	3.88406e-06\\
98525.44092	6.94e-07\\
116706.9647	4.67e-08\\
137016.9177	6.2e-08\\
153849.7427	6.53e-08\\
171217.5929	6.7e-08\\
};
\addlegendentry{Fine mesh}

\addplot [color=black, dashed, line width=1.pt] coordinates
{(6856.919206,0.000734635) (15708.29441,0.000512651) (54709.72303,0.000112548)};

\addplot [color=black, dashed, line width=1.pt] coordinates
{(8568.468921,0.000450835) (19982.28557,0.000208543) (68977.85379,2.28757e-05)};

\addplot [color=black, dashed, line width=1.pt] coordinates
{(10424.72281,0.000298845) (24311.20311,7.65471e-05) (83921.06135,3.88406e-06)};

\addplot [color=black, dashed, line width=1.pt] coordinates
{(12022.64435,0.000174209) (28541.96511,2.47091e-05) (98525.44092,6.94e-07)};

\addplot [color=black, dashed, line width=1.pt] coordinates
{(14497.4067,7.27122e-05) (33808.99476,7.01426e-06) (116706.9647,4.67e-08)};

\end{axis}

\end{tikzpicture}%
	\caption{Exponential convergence of Scordelis-Lo roof problem with three different meshes. Solid lines represent  errors generated from different polynomial orders (first to eighth) of the same mesh. Dashed lines represent errors from different meshes with the same polynomial order (from first to fifth).}
	\label{fig:scordelis_lo_roof_convergence}
\end{figure}
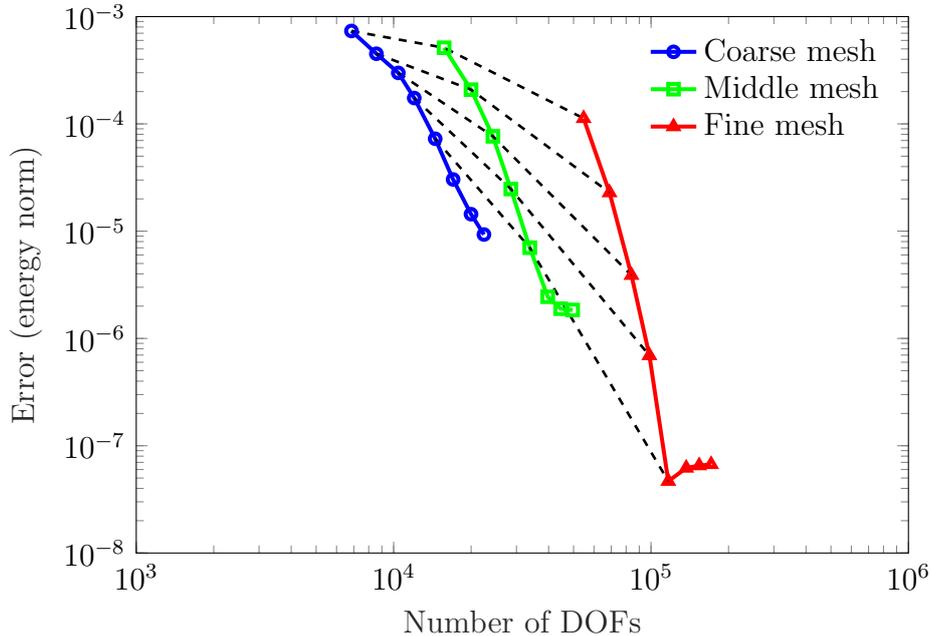

Fig. \ref{fig:scordelis_lo_roof_convergence} shows the convergence
of the analyses for Scordelis-Lo Roof with coarse, middle, and fine meshes for increasing order of approximation on the upper and lower faces. 
As can be seen, all three cases yield an exponential rate of convergence. The slight increase in error for the fine mesh for a very large number of DOFs
can be attributed to the tolerance of the linear solver, matrix conditioning, 
and the floating precision of the machine as the approximation error becomes small.

\begin{figure}[!htbp]
	\centering
	\includegraphics[scale=0.16]{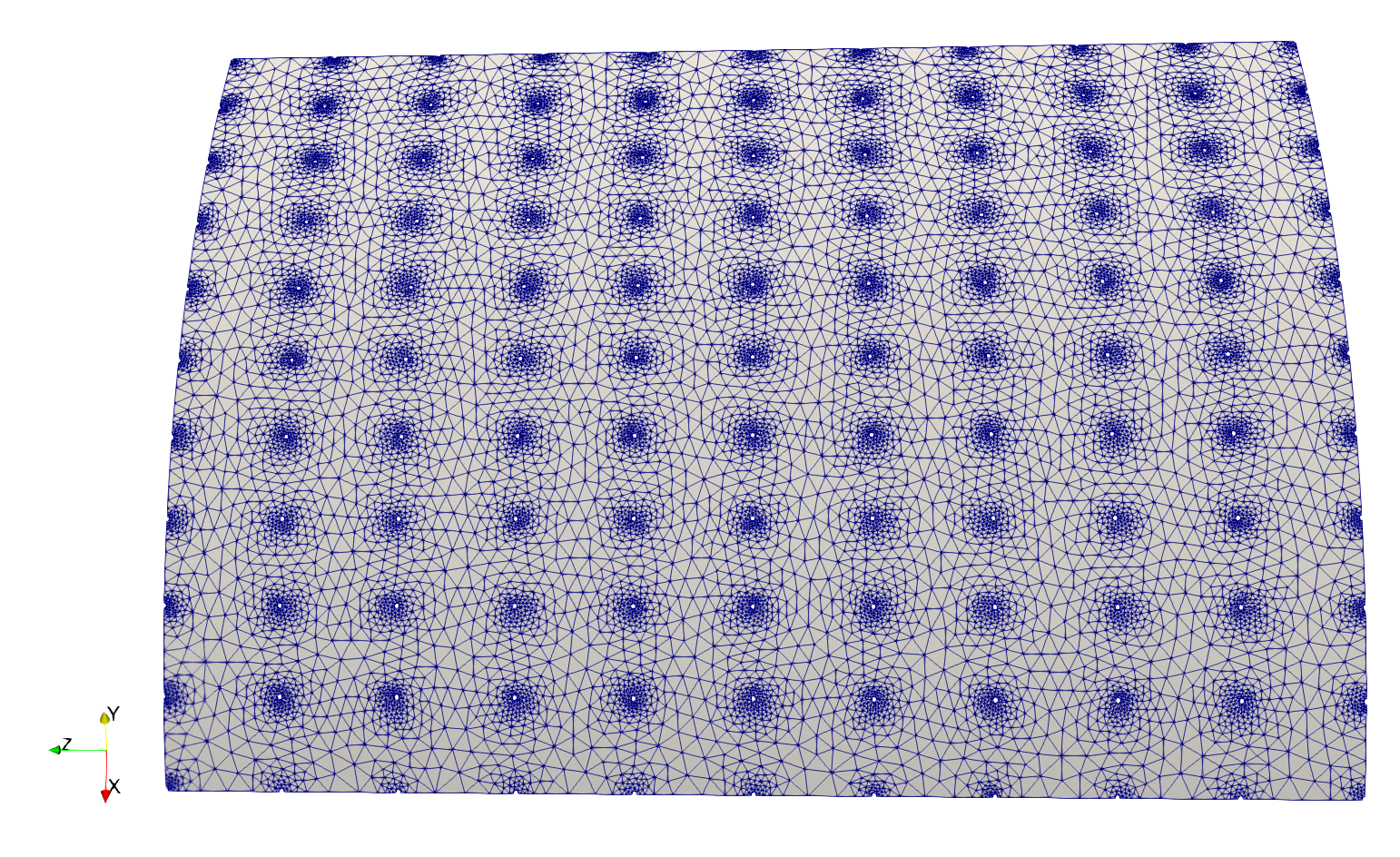}
	\caption{Mesh of a quarter of Scordelis-Lo roof with perforations.}
	\label{fig:scordelis_lo_roof_with_holes_geometry}
\end{figure}

\begin{figure}[!htbp]
	\begin{subfigure}{1.0\textwidth}
		\centering
		\includegraphics[scale=0.3]{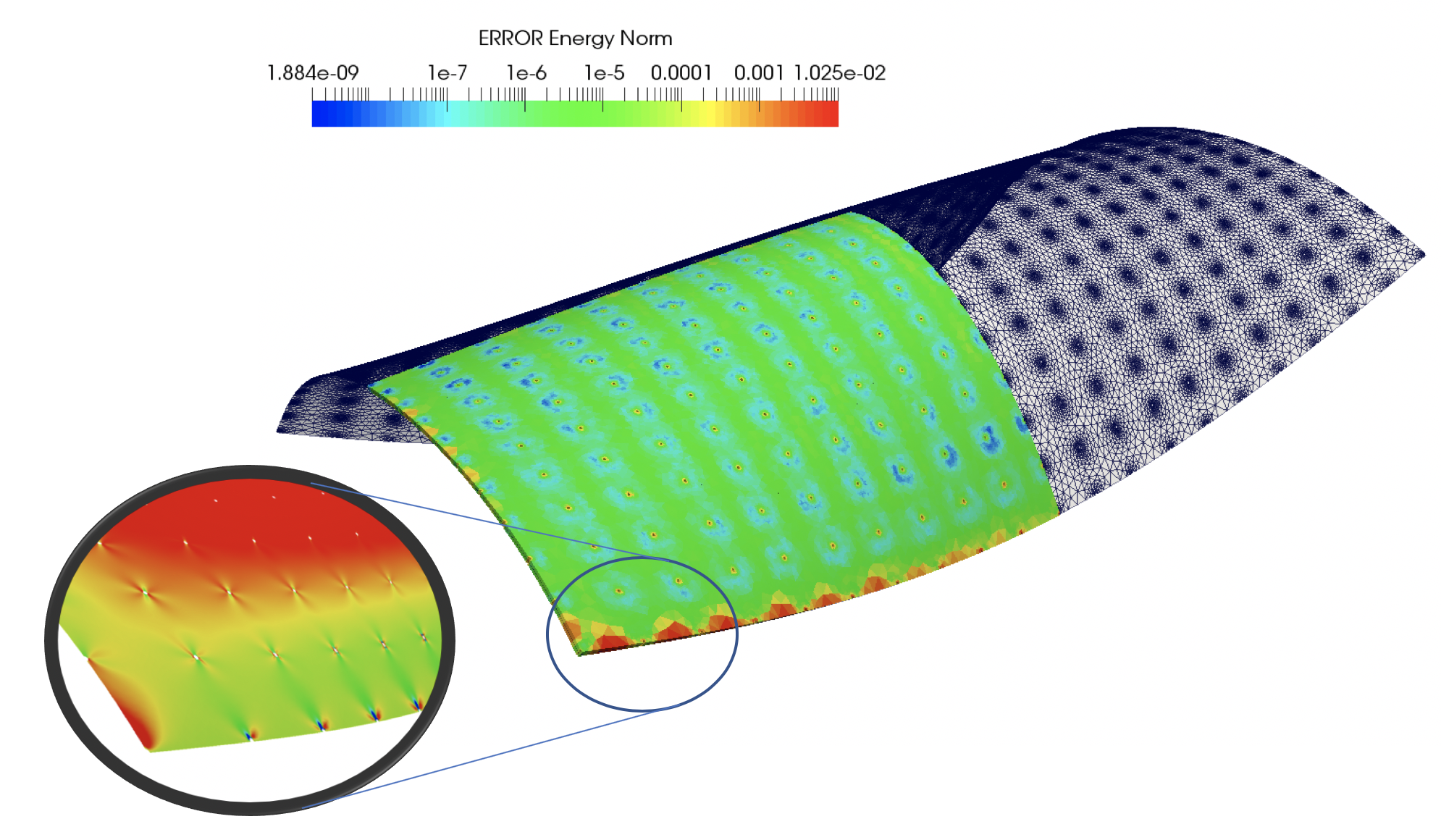}
		\caption{Error distribution before $p$-refinement}
		\label{fig:scordelis_lo_roof_with_holes_refinement_a}
	\end{subfigure}
	\begin{subfigure}{1.0\textwidth}
		\centering
		\includegraphics[scale=0.3]{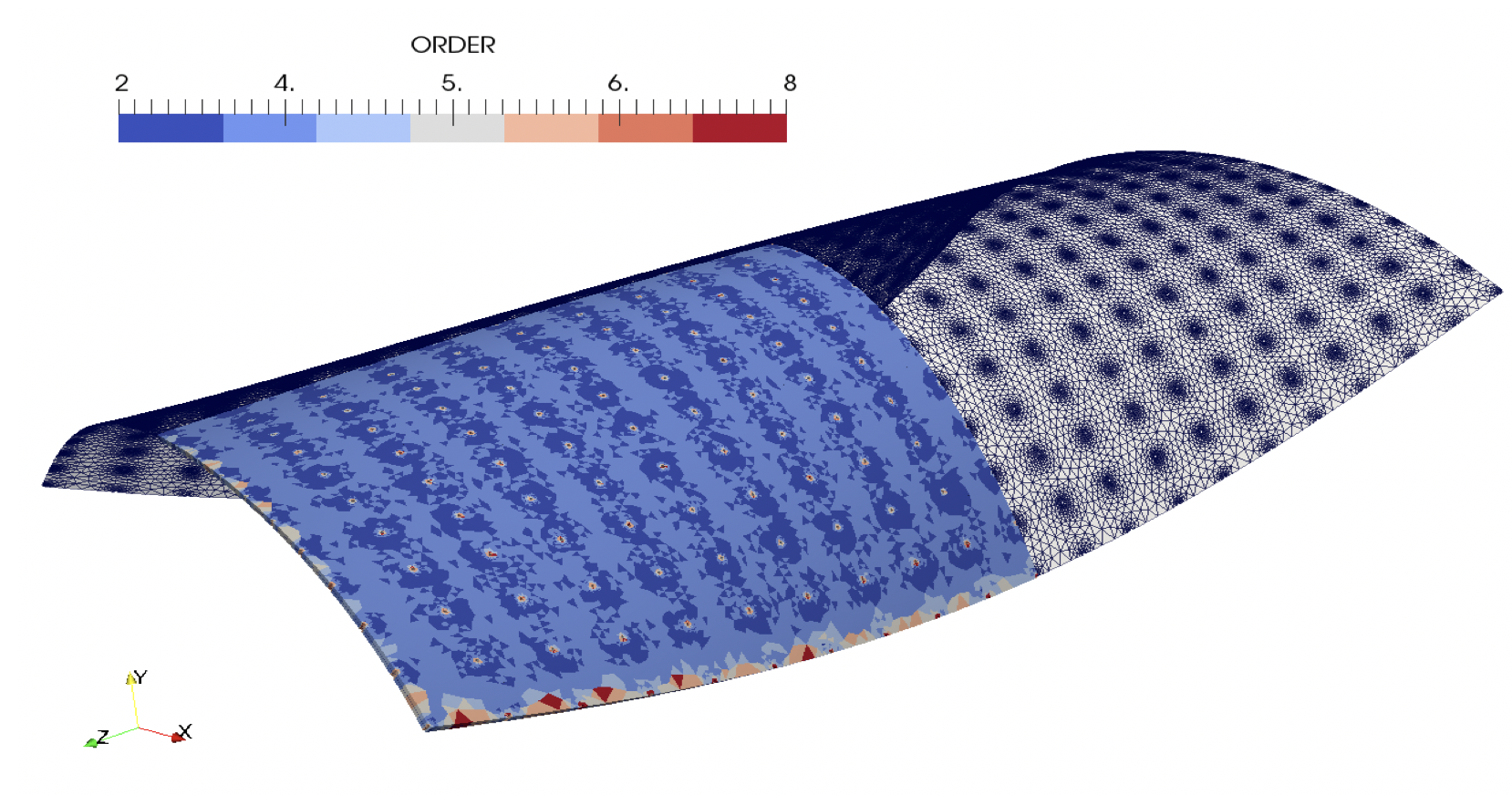}
		\caption{Heterogeneous local $p$-refinement with colour bar showing polynomial orders}
		\label{fig:scordelis_lo_roof_with_holes_refinement_b}
	\end{subfigure}
	\begin{subfigure}{1.0\textwidth}
		\centering
		\includegraphics[scale=0.3]{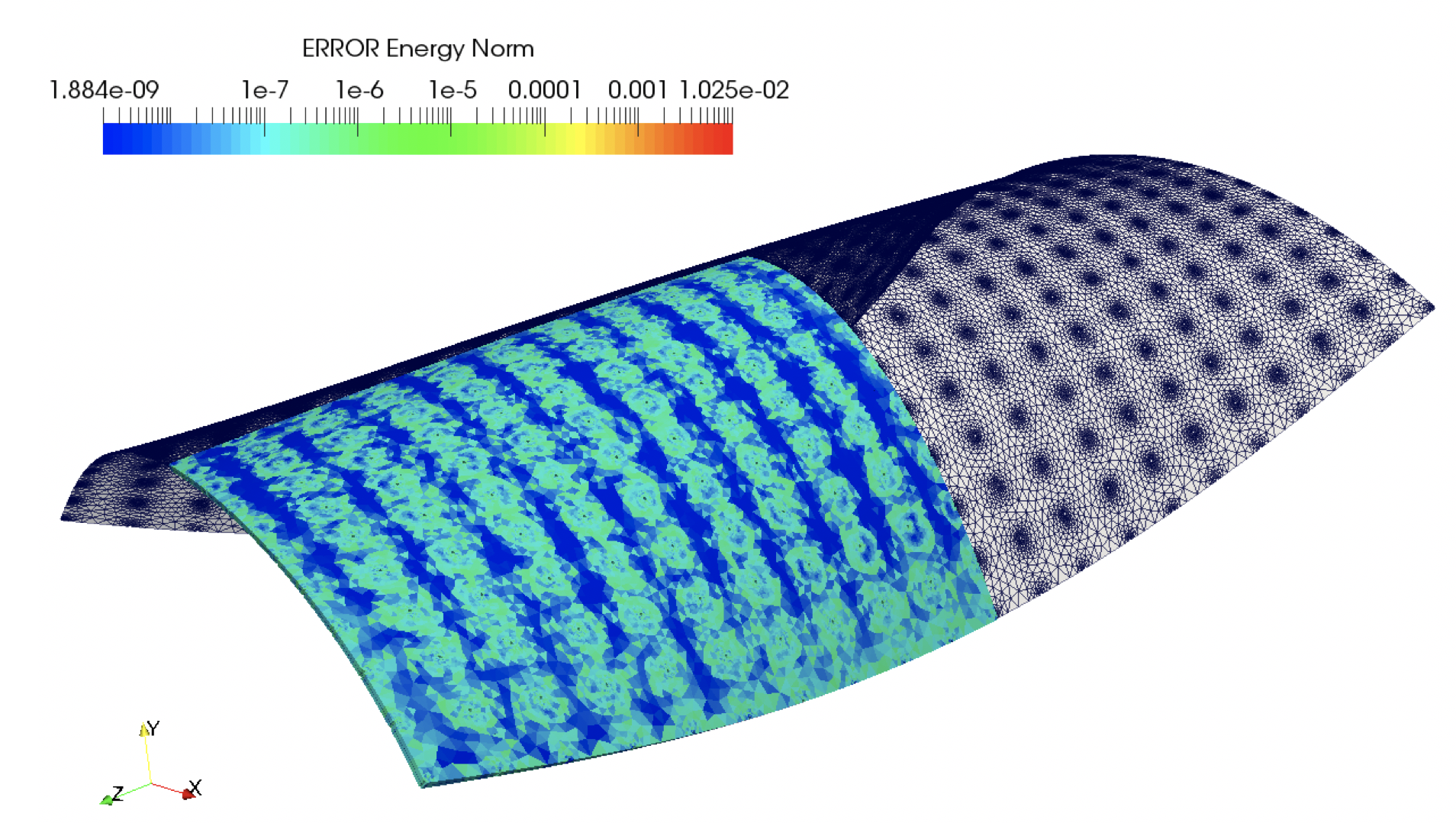}
		\caption{Error distribution after $p$-refinement}
		\label{fig:scordelis_lo_roof_with_holes_refinement_c}
  \end{subfigure}
	\caption{Multigrid solver with \emph{a posteriori} error estimator used for analysis
	of perforated Scordelis-Lo roof. 
		\label{fig:scordelis_lo_roof_with_holes_refinement}}
\end{figure}

The efficiency of the proposed prism elements based on hierarchical shape functions is demonstrated using a more
complex geometry of the Scordelis-Lo Roof, with perforations. Fig.
\ref{fig:scordelis_lo_roof_with_holes_geometry} shows this modified geometry
with uniformly distributed holes. Each hole has a radius of $0.3$.
Once again, due to symmetry, only a quarter of the perforated roof is considered. 
The perforated roof is initially
analysed using second-order approximations on the triangular faces, 
yielding poor results, as expected. The error in energy norm is 
distributed non-uniformly across the roof with the largest errors 
observed near the outer edge
of the structure as shown in Fig.
\ref{fig:scordelis_lo_roof_with_holes_refinement}. 
This \emph{a posteriori} error estimator is subsequently employed to 
drive adaptivity, in which one third of the elements with the highest 
error are selected to locally increase the order of approximation. 
This creates a heterogeneous distribution of approximation, 
ranging from second to eighth-order. 
An iterative adaptive refinement process is undertaken, leading to a highly
accurate result with low errors across the roof.

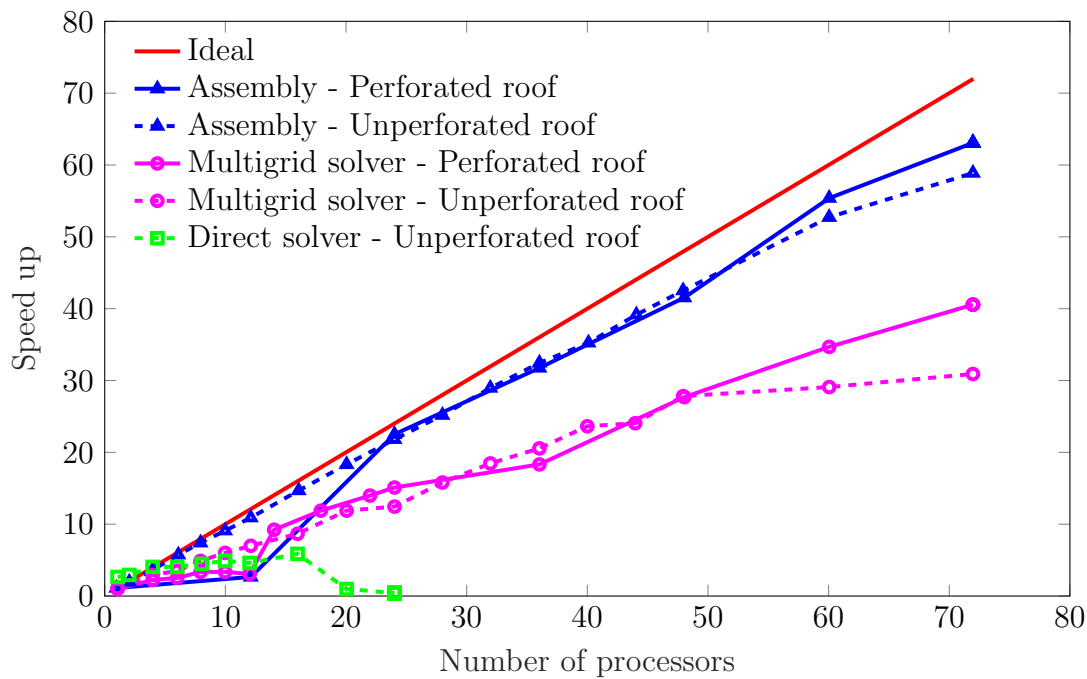
\begin{figure}[!htb]
	\centering
	~
%
%
\definecolor{mycolor1}{rgb}{1.00000,0.00000,1.00000}%
\begin{tikzpicture}

\begin{axis}[%
width=5.0in,
height=3.0in,
scale only axis,
xmin=0,
xmax=80,
xtick={0,10,20,30,40,50,60,70,80},
xlabel style={font=\color{white!15!black}},
xlabel={Number of processors},
ymin=0,
ymax=80,
ylabel style={font=\color{white!15!black}},
ylabel={Speed up},
axis background/.style={fill=white},
legend style={at={(0.02,0.58)}, anchor=south west, legend cell align=left, align=left, fill=none, 
draw=none}
]
\addplot [color=red, line width=1.5pt]
  table[row sep=crcr]{%
1	1\\
72	72\\
72	72\\
72	72\\
72	72\\
72	72\\
72	72\\
72	72\\
72	72\\
72	72\\
72	72\\
72	72\\
72	72\\
72	72\\
72	72\\
72	72\\
72	72\\
};
\addlegendentry{Ideal}

\addplot [color=blue, line width=1.5pt, mark=triangle, mark options={solid, blue}]
  table[row sep=crcr]{%
0.968523002	1.118881119\\
12.10653753	2.657342657\\
24.01937046	22.51748252\\
36.02905569	31.74825175\\
48.03874092	41.53846154\\
60.04842615	55.38461538\\
71.96125908	63.07692308\\
71.96125908	63.07692308\\
71.96125908	63.07692308\\
71.96125908	63.07692308\\
71.96125908	63.07692308\\
71.96125908	63.07692308\\
71.96125908	63.07692308\\
71.96125908	63.07692308\\
71.96125908	63.07692308\\
71.96125908	63.07692308\\
71.96125908	63.07692308\\
};
\addlegendentry{Assembly - Perforated roof}

\addplot [color=blue, dashed, line width=1.5pt, mark=triangle, mark options={solid, blue}]
  table[row sep=crcr]{%
2.033898305	1.958041958\\
3.97094431	3.636363636\\
6.101694915	5.734265734\\
7.94188862	7.412587413\\
9.975786925	9.090909091\\
12.10653753	10.90909091\\
16.07748184	14.68531469\\
20.04842615	18.32167832\\
24.01937046	21.81818182\\
27.99031477	25.17482517\\
31.96125908	28.95104895\\
36.02905569	32.44755245\\
40.0968523	35.24475524\\
44.06779661	39.16083916\\
47.94188862	42.51748252\\
60.04842615	52.72727273\\
71.96125908	58.88111888\\
};
\addlegendentry{Assembly - Unperforated roof}

\addplot [color=mycolor1, line width=1.5pt, mark=o, mark options={solid, mycolor1}]
  table[row sep=crcr]{%
1.065375303	0.979020979\\
3.97094431	2.237762238\\
6.004842615	2.517482517\\
7.94188862	3.356643357\\
9.975786925	3.356643357\\
12.00968523	3.076923077\\
14.04358354	9.230769231\\
17.91767554	11.88811189\\
21.98547215	13.98601399\\
24.01937046	15.1048951\\
36.02905569	18.32167832\\
48.03874092	27.69230769\\
60.04842615	34.68531469\\
71.96125908	40.55944056\\
71.96125908	40.55944056\\
71.96125908	40.55944056\\
71.96125908	40.55944056\\
};
\addlegendentry{Multigrid solver - Perforated roof}

\addplot [color=mycolor1, dashed, line width=1.5pt, mark=o, mark options={solid, mycolor1}]
  table[row sep=crcr]{%
1.065375303	1.118881119\\
3.97094431	2.937062937\\
6.004842615	3.496503497\\
7.94188862	4.895104895\\
9.975786925	6.013986014\\
12.10653753	6.993006993\\
15.98062954	8.671328671\\
20.04842615	11.88811189\\
24.01937046	12.44755245\\
27.99031477	15.8041958\\
31.96125908	18.46153846\\
36.02905569	20.55944056\\
40	23.63636364\\
43.97094431	24.05594406\\
47.94188862	27.83216783\\
60.04842615	29.09090909\\
71.96125908	30.90909091\\
};
\addlegendentry{Multigrid solver - Unperforated roof}

\addplot [color=green, dashed, line width=1.5pt, mark=square, mark options={solid, green}]
  table[row sep=crcr]{%
1.065375303	2.657342657\\
2.033898305	2.937062937\\
3.97094431	4.055944056\\
6.004842615	4.055944056\\
8.03874092	4.475524476\\
9.975786925	4.895104895\\
12.00968523	4.615384615\\
15.98062954	5.874125874\\
20.04842615	0.979020979\\
24.01937046	0.41958042\\
24.01937046	0.41958042\\
24.01937046	0.41958042\\
24.01937046	0.41958042\\
24.01937046	0.41958042\\
24.01937046	0.41958042\\
24.01937046	0.41958042\\
24.01937046	0.41958042\\
};
\addlegendentry{Direct solver - Unperforated roof}

\end{axis}

\begin{axis}[%
width=7.778in,
height=5.833in,
at={(0in,0in)},
scale only axis,
xmin=0,
xmax=1,
ymin=0,
ymax=1,
axis line style={draw=none},
ticks=none,
axis x line*=bottom,
axis y line*=left,
legend style={legend cell align=left, align=left, draw=white!15!black}
]
\end{axis}
\end{tikzpicture}%
	\caption{Speedup for assembly and solver for unperforated and perforated Scordelis-Lo roof.}
	\label{fig:scordelis_lo_roof_with_holes_multigrid_scalability}
\end{figure}

It is worth noting that, in this example of linear analysis, that a multigrid
solver is utilised to significantly 
increase the speed of solving the system of equations. Due to its scalability,
the multigrid solver, along with hierarchical shape functions and an error estimator,
appears to be a highly effective combination in achieving an efficient analysis
using multiple processors. Furthermore, the anisotropic
and heterogeneous properties of the shape functions support more
flexibility, leading to accurate results with minimum additional DOFs. Fig.
\ref{fig:scordelis_lo_roof_with_holes_multigrid_scalability} shows that the 
multigrid solver scales well with increasing numbers of processors. It also shows the excellent scalability of the assembly process.

For the unperforated roof, due to the relative simplicity of the problem 
and small number of DOFs, there is a sign of
saturation. In contrast, a direct solver (SuperLU\_DIST \cite{li_overview_2005})
performs slightly better for a small number of processors. However, its performance quickly deteriorates as the number of processors grow. This can be attributed to increased matrix fill-in when the problem is distributed over many processors. After a certain point, with more processors, the speed-up for the direct solver even decreases.

\subsection{Geometrically nonlinear analysis}

\subsubsection{Slit annular plate}
\label{sec:annular_plate}
\begin{figure}[!htbp]
	\centering
	\def\svgwidth{8cm}
	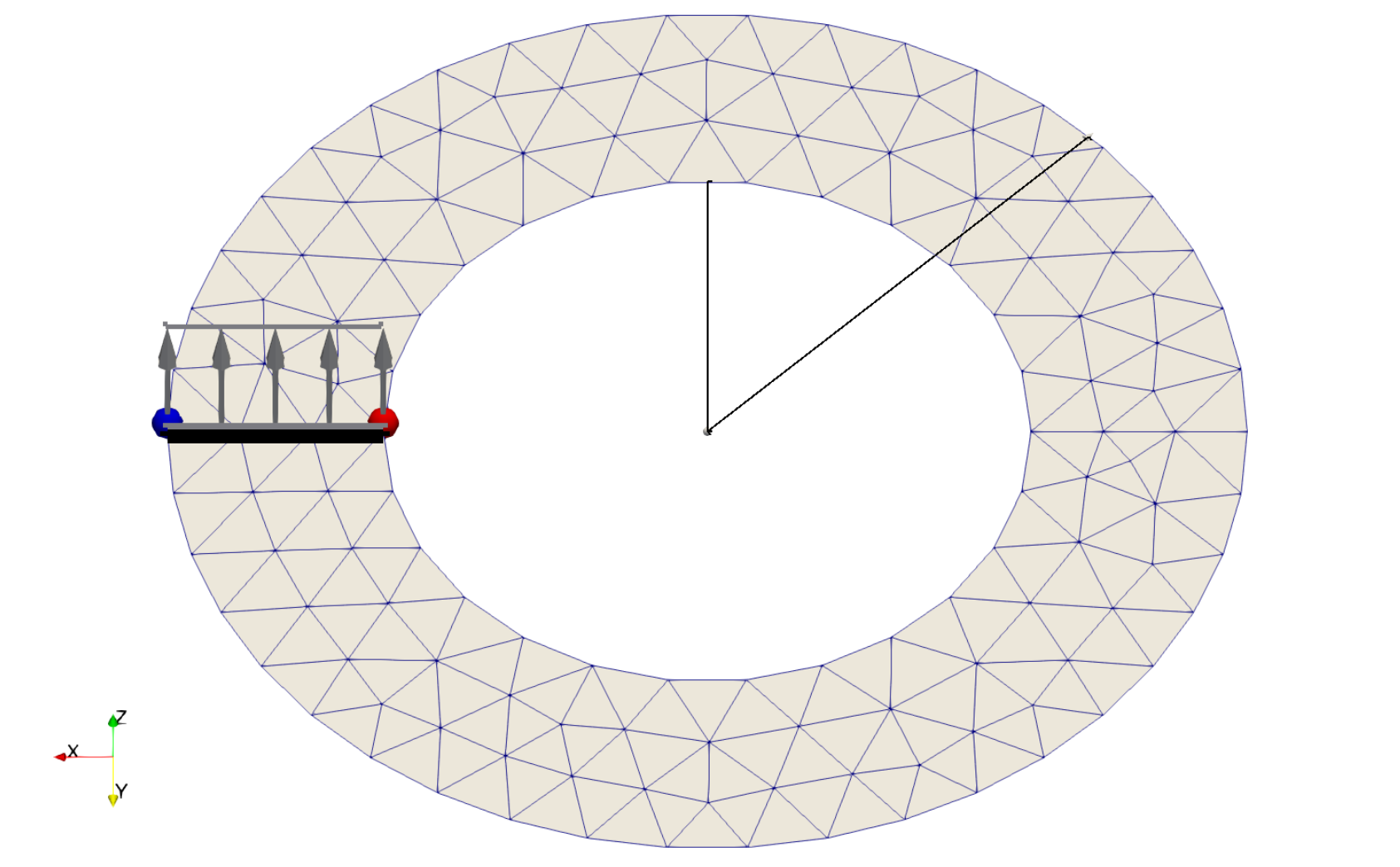	
	\caption{Configuration of slit annular plate.}
	\label{fig:slit_annular_plate_geometry}
\end{figure}

\begin{figure}[!htbp]
	\centering
	~\input{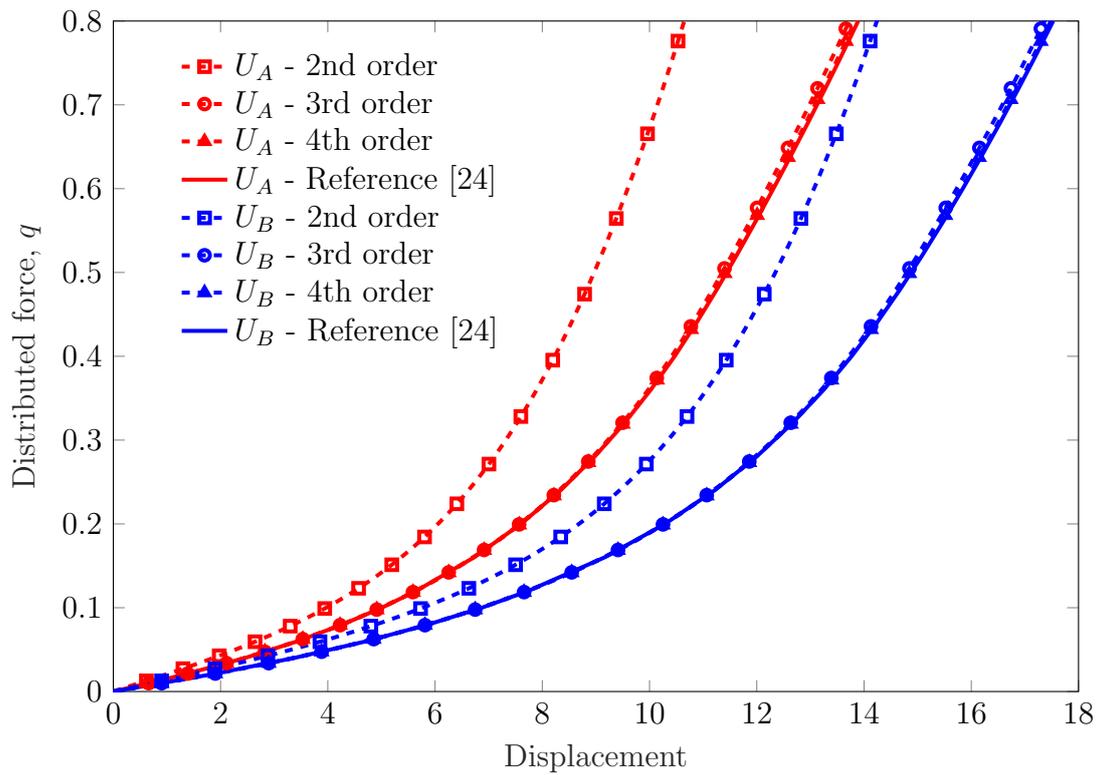}
	\caption{Comparison of slit annular plate at observed points with different orders of approximation for the upper and lower faces.}
	\label{fig:slit_annular_plate_comparison}
\end{figure}

The nonlinear analysis of a slit annular plate, subject to an edge
force, is considered as a standard benchmark for the nonlinear analysis shells. Fig. \ref{fig:slit_annular_plate_geometry} shows the
geometry and mesh for this problem. The input data is
provided in Table \ref{tab:input_data}. One side of the slit is fully
clamped and the other side is subjected a uniformly distributed edge load which is
kept vertical during the entire simulation. The vertical
displacements of point A (red) and point B (blue) at two ends of the edge
are observed and compared to those in the literature.
Fig. \ref{fig:slit_annular_plate_comparison} presents the displacements at both points 
for increasing force, comparing different order of approximation on the upper and lower faces. 
Once again, second-order approximations are adopted through the thickness of the element.

\begin{figure}[!htb]
	\begin{subfigure}{1.0\textwidth}
		\centering
		\includegraphics[scale=0.4]{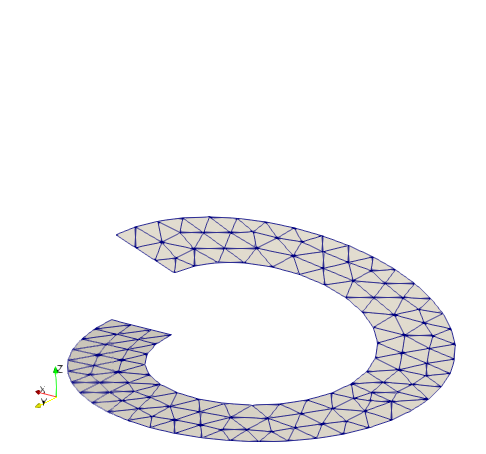}
		\caption{Step 50 ($q=0.064$)}
	\end{subfigure}
	\begin{subfigure}{1.0\textwidth}
		\centering
		\includegraphics[scale=0.4]{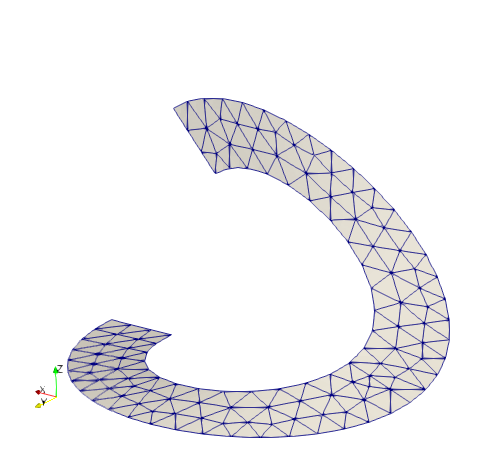}
		\caption{Step 150 ($q=0.380$)}
	\end{subfigure}
	\begin{subfigure}{1.0\textwidth}
		\centering
		\includegraphics[scale=0.4]{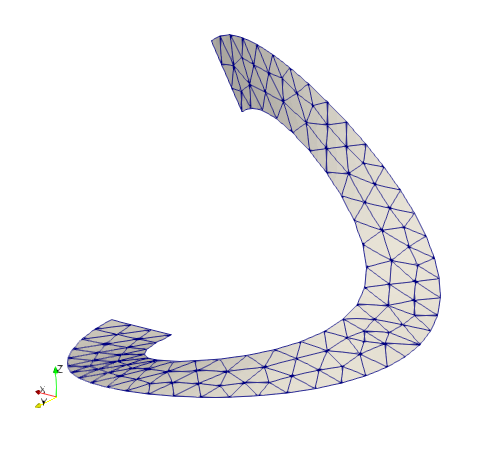}
		\caption{Step 211 ($q=0.8$)}
	\end{subfigure}
	\caption{Deformed shapes for the nonlinear analysis of slit annular plate at different stages of the analysis.
		\label{fig:slit_annular_plate_steps}}
\end{figure}

It can be observed that a second-order approximation on the top and bottom triangles yields substantially stiffer results compared to the results reported by Sze et al.
\cite{sze_popular_2004}. This volume locking can be attributed to the coarse mesh and low order of approximation employed. 
However, as the approximation is
raised to third and fourth orders, locking is alleviated and the
vertical displacements at points A and B agree very well with the reference. The
deformed shapes of the plate for third order approximation are
shown in Fig. \ref{fig:slit_annular_plate_steps}.

\subsubsection{Hemispherical shell}
\label{sec:hemispherical_shell}

\begin{figure}[!htbp]
	\centering
	\def\svgwidth{8cm}
	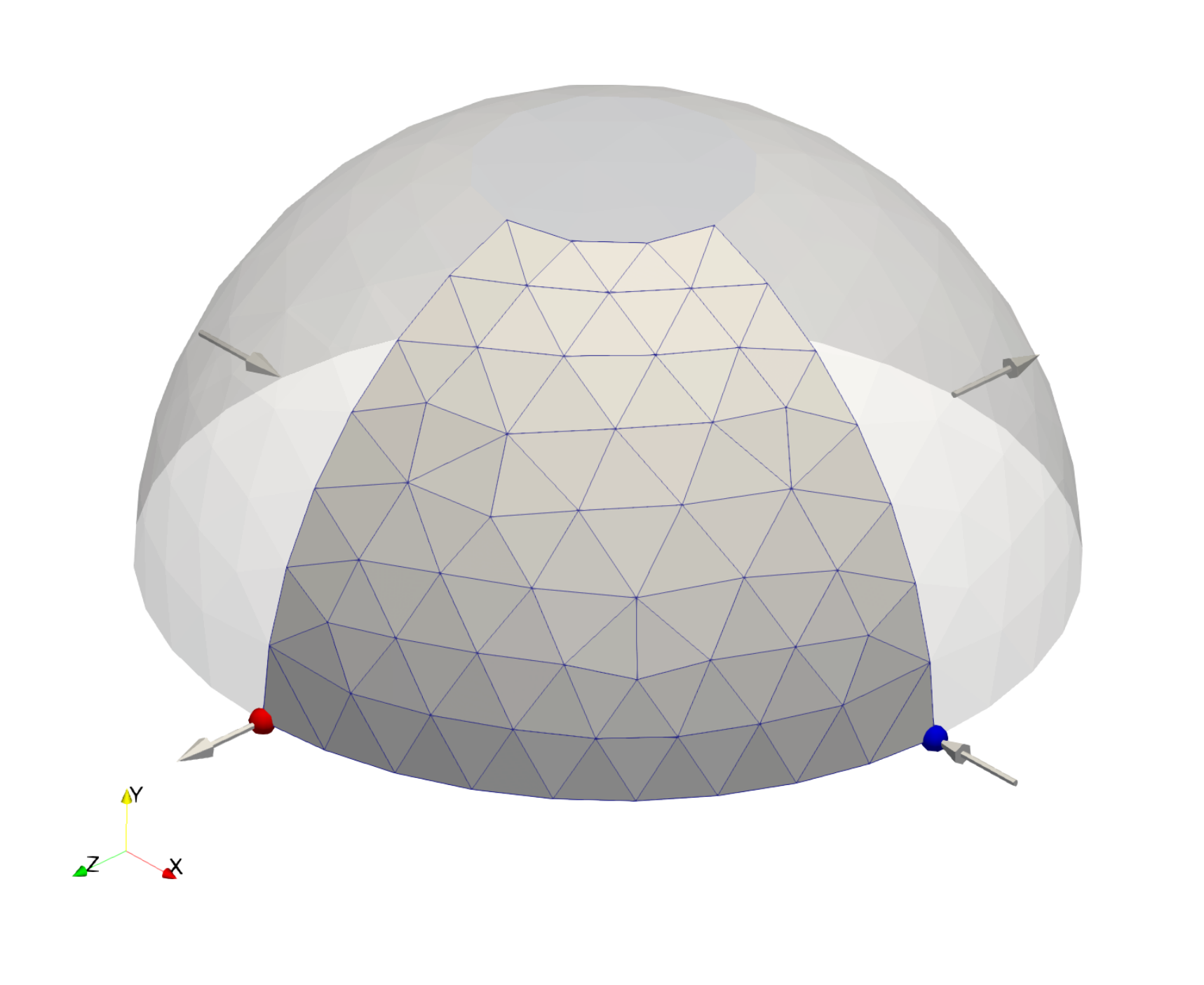	
	\caption{Geometry, mesh and loading of hemispherical shell}
	\label{fig:hemispherical_shell_geometry}
\end{figure}

Fig. \ref{fig:hemispherical_shell_geometry} shows the geometry and mesh
configuration of a hemispherical shell with an $18^\circ$ cutout at its pole. 
Detailed dimensions and material properties of the structure can
be found in Table \ref{tab:input_data}. Due to symmetry, only one quarter of
the shell with symmetric boundary conditions is analysed. Vertical displacement is
restricted along the bottom edge.
The radial displacement at points A and B are observed and compared with 
the reference \cite{sze_popular_2004}. This hemispherical shell problem is 
considered a valuable test challenging test due to the inextensional 
bending modes with almost no membrane strains.
Further, since a large portion of
the shell rotates under the given loads, it is also a useful test to evaluate the capability of the element in handling rigid body
rotation about normals to the shell surfaces \cite{belytschko_stress_1985}. 

\begin{figure}[!htbp]
	\centering
	~\input{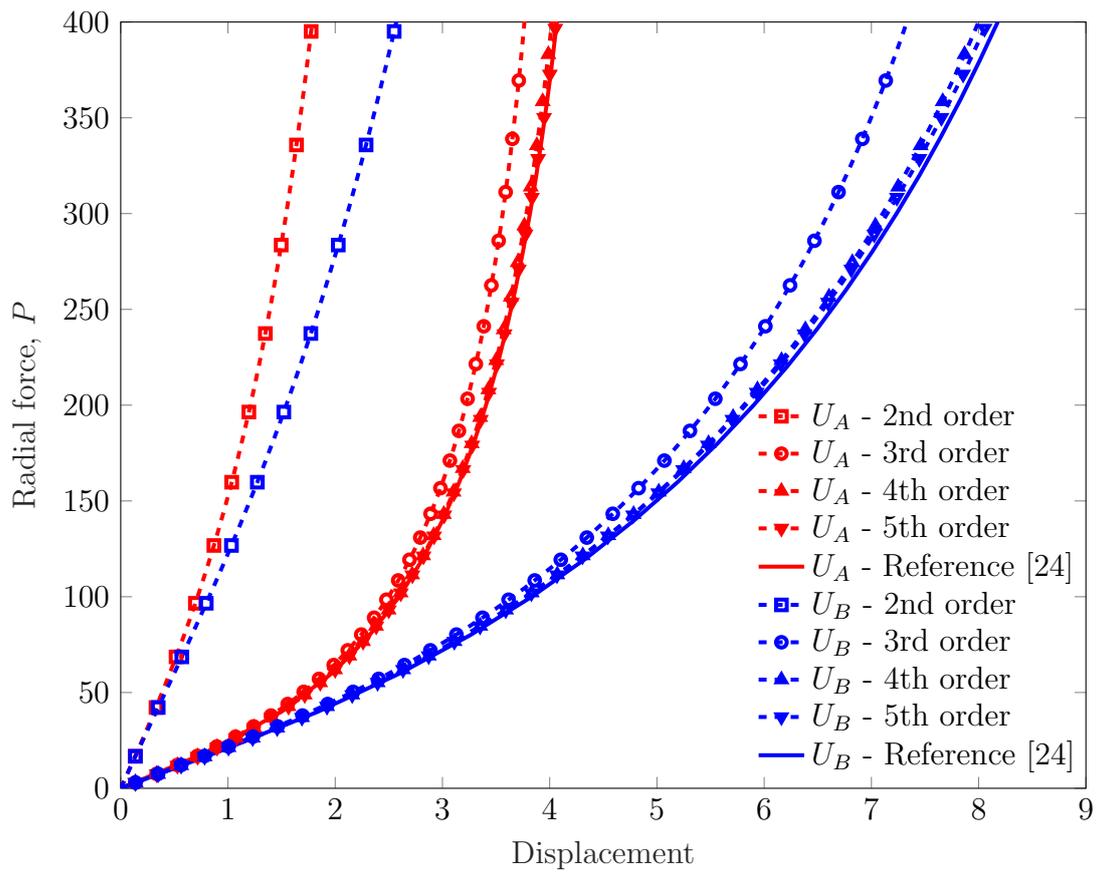}
	\caption{Comparison of radial displacements of hemispherical shell at points A (red) and B (blue) for different orders of approximation.}
	\label{fig:hemispherical_shell_comparison}
\end{figure}

\begin{figure}[!htbp]
	\begin{subfigure}{1.0\textwidth}
		\centering
		\includegraphics[scale=0.3]{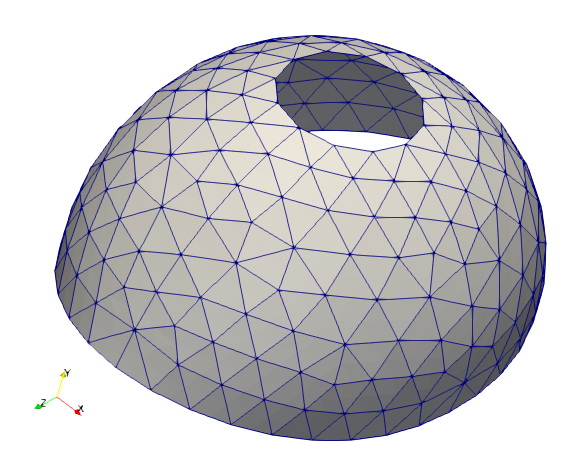}
		\caption{Step 20 ($P=32$)}
	\end{subfigure}
	\begin{subfigure}{1.0\textwidth}
		\centering
		\includegraphics[scale=0.3]{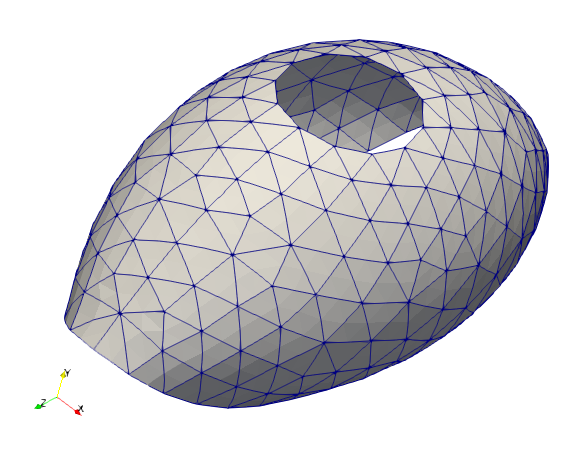}
		\caption{Step 70 ($P=198$)}
	\end{subfigure}
	\begin{subfigure}{1.0\textwidth}
		\centering
		\includegraphics[scale=0.3]{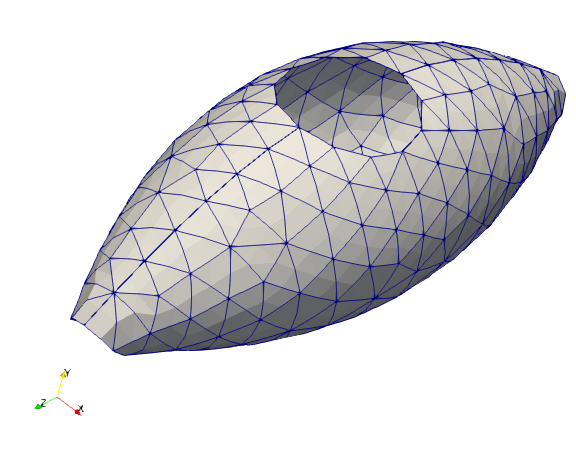}
		\caption{Step 137 (last step, $P=1,578$)}
	\end{subfigure}
	\caption{Deformed shapes of hemispherical shell at three different stages of the analysis.
		\label{fig:hemispherical_shell_steps}}
\end{figure}

Similar to the previous example, a coarse mesh with only 100 elements is used to
analyse the shell and evaluate the performance of the proposed approach, with second-order
approximations through the thickness.
The approximations constructed on the upper and lower faces are varied from second to fifth orders.
As in the previous problem, a lower order of approximation leads to
volumetric locking and a stiff result, as shown in Fig.
\ref{fig:hemispherical_shell_comparison}. However, the locking appears to be 
alleviated when
fourth and fifth orders of approximation are employed and the results are in good
agreement with those reported in the literature \cite{sze_popular_2004} for both
points A and B. 
Fig.
\ref{fig:hemispherical_shell_steps} presents the deformed 
shapes for third-order approximations.

\subsubsection{Cylindrical shell bearing pull-out forces}
\label{sec:cylindrical_open_ended}

\begin{figure}[!htbp]
	\centering
	\def\svgwidth{9cm}
	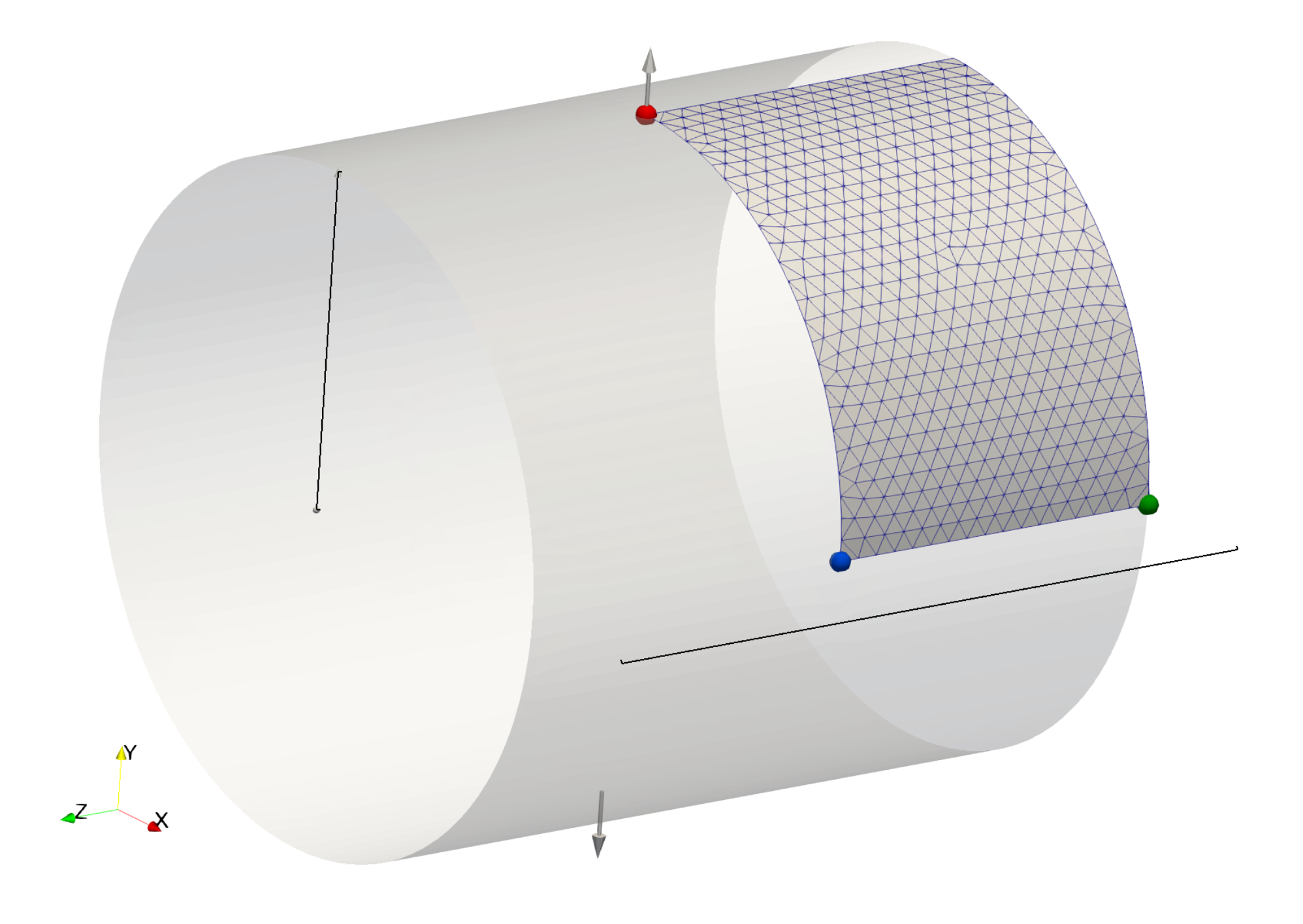	
	\caption{Geometry, boundary conditions, loading and mesh for the pull-out cylinder.}
	\label{fig:pullout_cylinder_geometry}
\end{figure}

Fig. \ref{fig:pullout_cylinder_geometry} shows the geometry and mesh of an
open-ended cylindrical shell subjected to a pair of radial forces $P$. Due to 
symmetry, only one eighth of the shell is analysed and the radial
displacement of points A, B, and C are observed and compared to the reference
\cite{sze_popular_2004}. The geometry and material data can be found
in Table \ref{tab:input_data}.

\begin{figure}[!htbp]
	\centering
	~\input{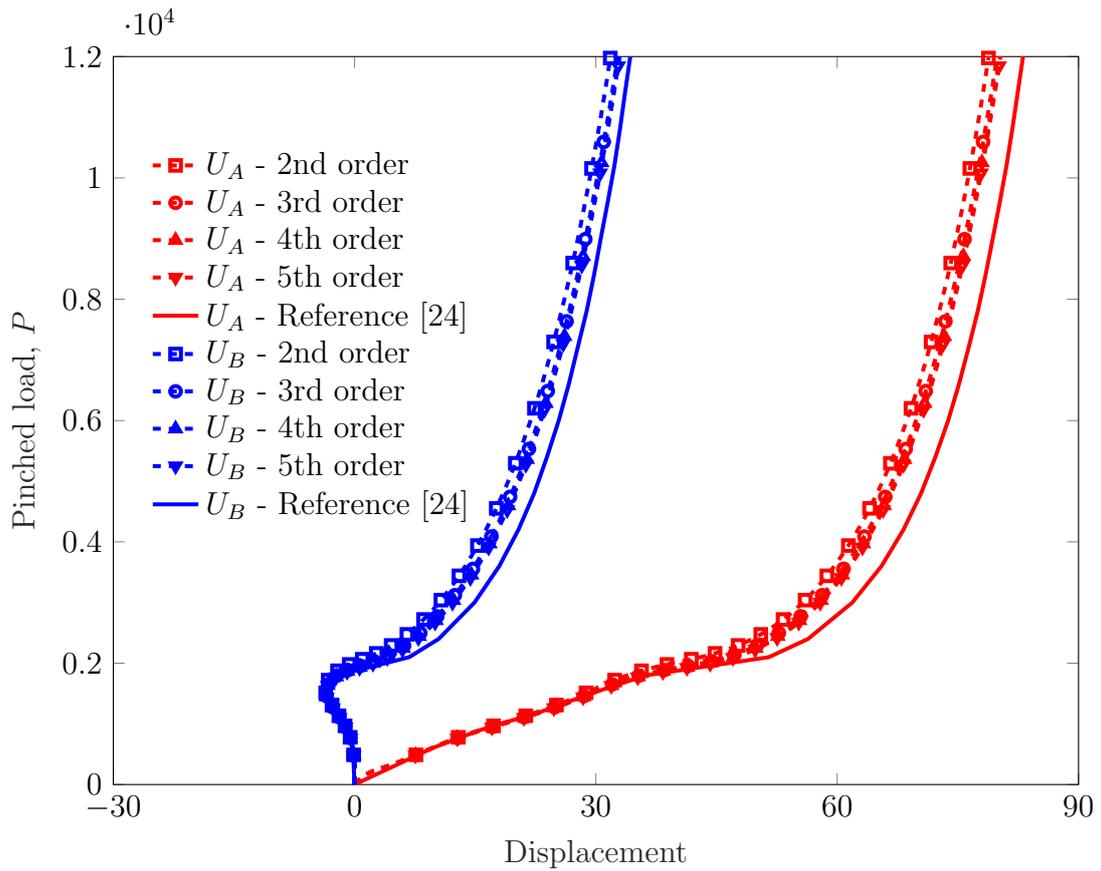}
	\caption{Comparison of radial displacements of the pull-out cylindrical
	shell at points A (red), B (blue), and C (green) for different orders of approximation.
	}
	\label{fig:pullout_cylinder_comparison}
\end{figure}

Due to the relatively fine mesh being used, the radial displacement of the 
observed points
shows good agreement with those reported by Sze et al. \cite{sze_popular_2004},
even with second or third orders of approximation on the prism upper and lower faces. 
This is illustrated in load-displacement plot in Fig \ref{fig:pullout_cylinder_comparison}. 

\begin{figure}[!htb]
	\begin{subfigure}{1.0\textwidth}
		\centering
		\includegraphics[scale=0.4]{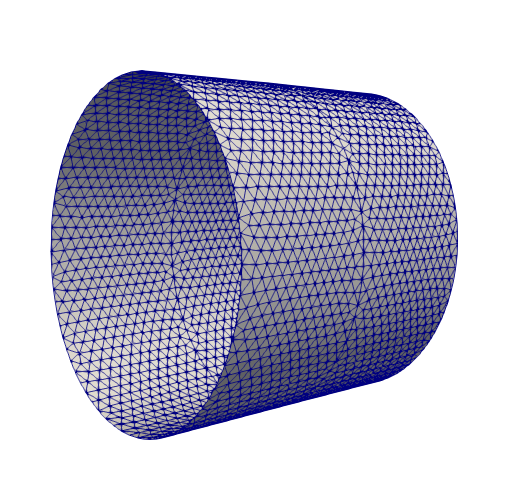}
		\caption{Step 10 ($P=528$)}
	\end{subfigure}
	\begin{subfigure}{1.0\textwidth}
		\centering
		\includegraphics[scale=0.4]{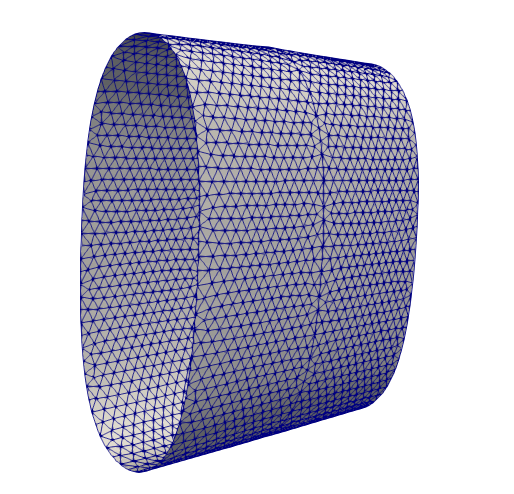}
		\caption{Step 40 ($P=3,671$)}
	\end{subfigure}
	\begin{subfigure}{1.0\textwidth}
		\centering
		\includegraphics[scale=0.4]{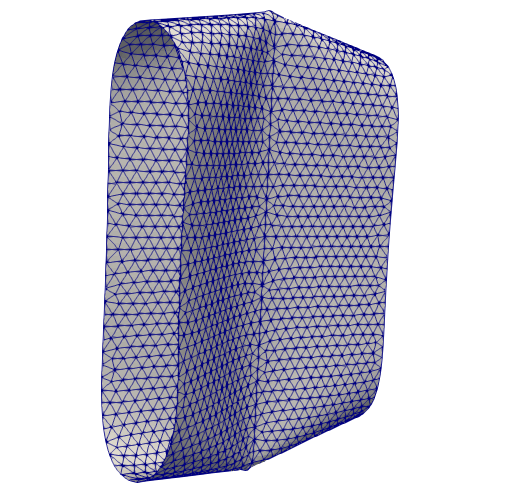}
		\caption{Step 87 (last step, $P=52,414$)}
	\end{subfigure}
	\caption{Deformed shapes of the cylindrical shell at three different stages of the analysis. 
		\label{fig:pullout_shell_steps}}
\end{figure}

The large and complex displacements of point C at the edge of the shell
is captured well by the proposed approach.
The deformed shapes at three analysis steps with third-order approximations are presented in Fig. \ref{fig:pullout_shell_steps}. 
Second-order approximations were adopted in the through-thickness direction for all analyses.


\subsubsection{Pinched cylindrical shell}
\label{sec:cylindrical_rigid_diaphragm}

\begin{figure}[!htbp]
	\centering
	\def\svgwidth{9cm}
	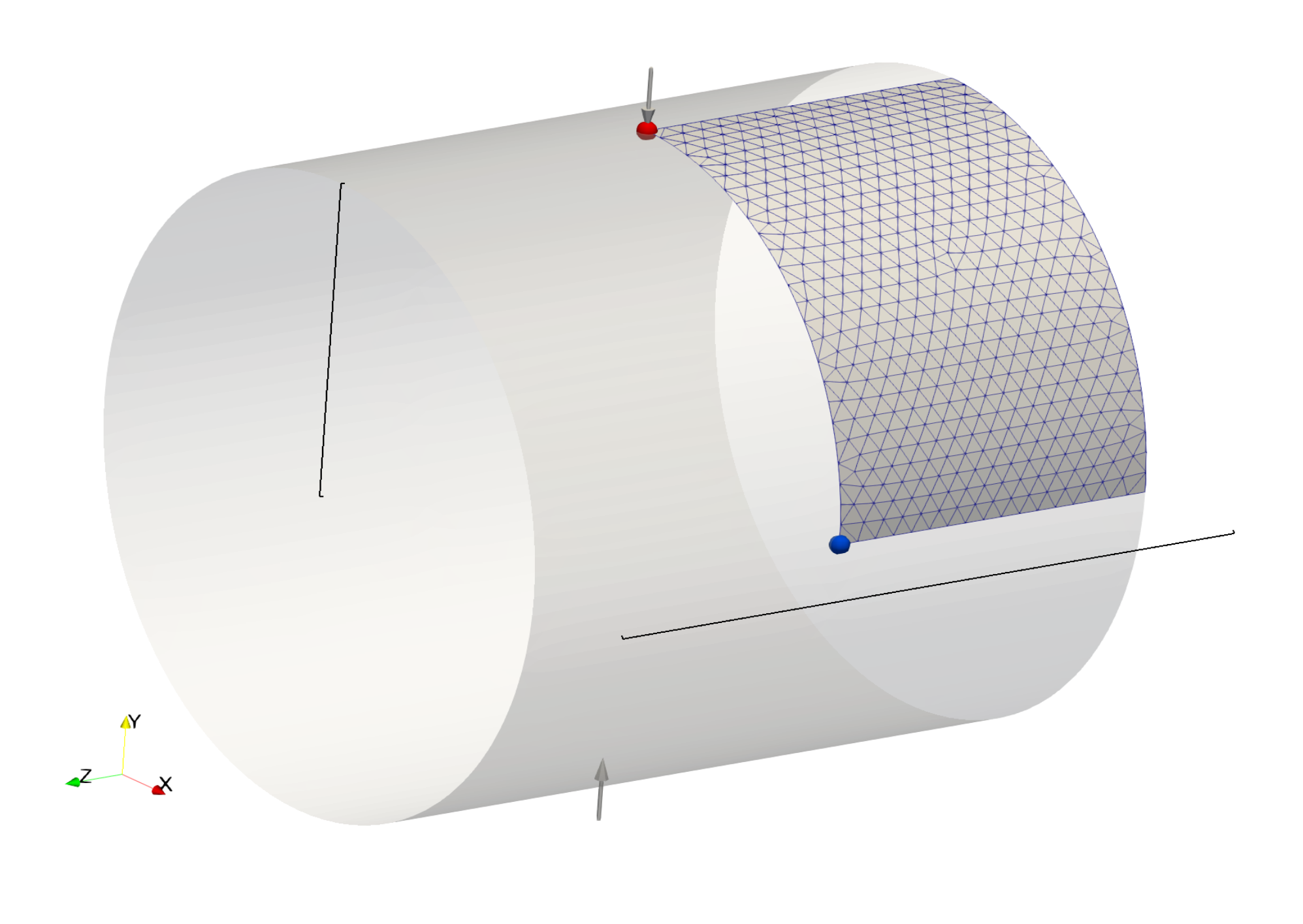	
	\caption{Geometry, boundary conditions, loading and mesh for pinched 
	cylindrical shell with rigid diaphragms.}
	\label{fig:pinched_cylinder_geometry}
\end{figure}

This example is similar to the previous but the cylindrical shell is
mounted on to rigid diaphragms at each end, restricting in-plane displacements. The
geometry and mesh configuration and its details are  presented in Fig.
\ref{fig:pinched_cylinder_geometry} and Table \ref{tab:input_data},respectively.
Once again, due to symmetry, one eighth of the cylinder is analysed. This test is 
particularly challenging due to the need to capture 
both the inextensional bending modes and the complex membrane states 
which the
cylinder exhibits during deformation \cite{belytschko_stress_1985}.

\begin{figure}[!htbp]
	\centering
	~\input{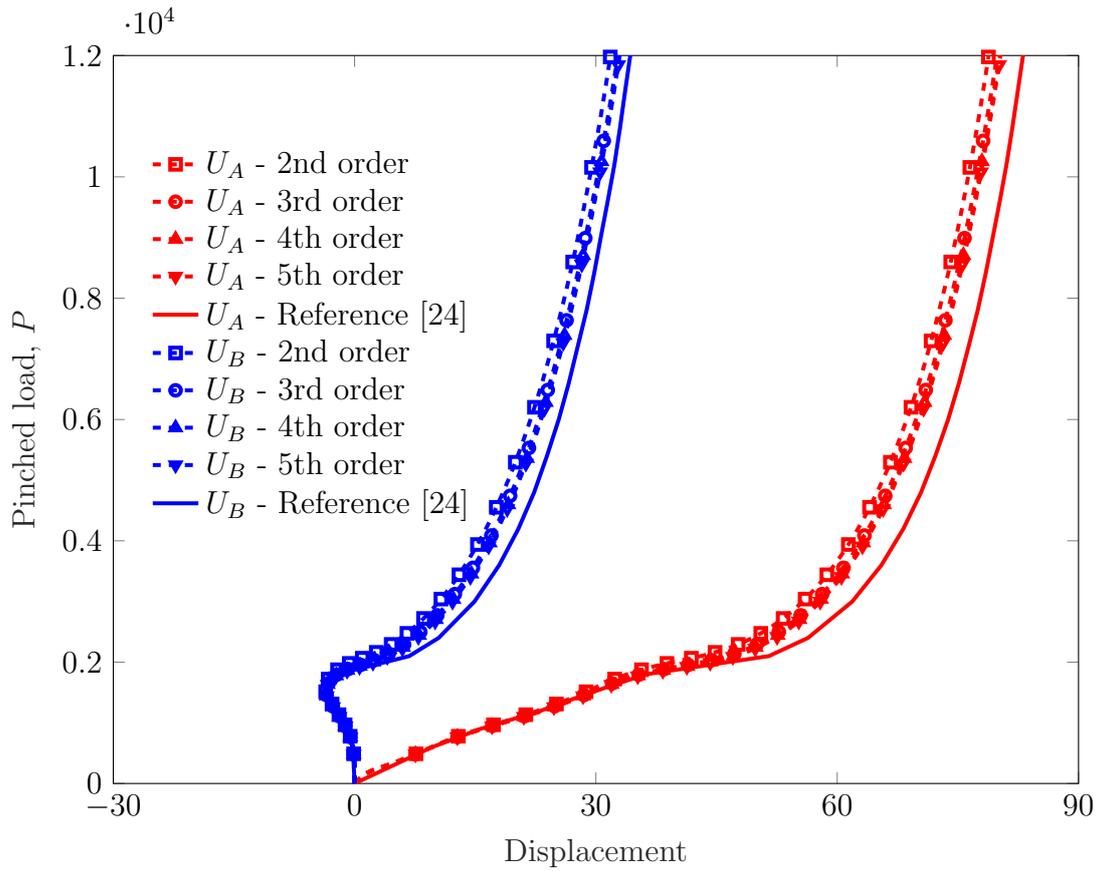}
	\caption{Comparison of radial displacements of the pinched shell at points A (red) and B (blue) for different orders of approximation.}
	\label{fig:pinched_cylinder_comparison}
\end{figure}

\begin{figure}[!htbp]
	\begin{subfigure}{1.0\textwidth}
		\centering
		\includegraphics[scale=0.3]{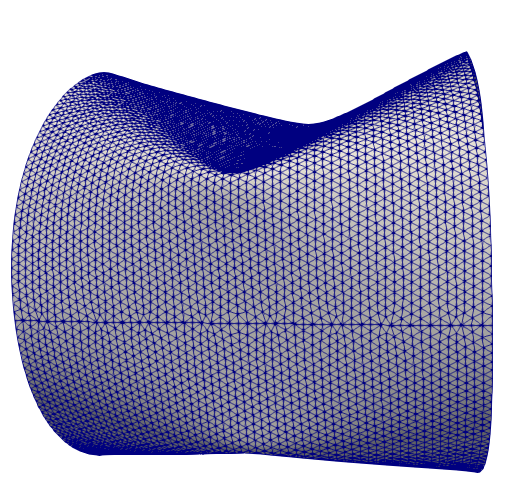}
		\caption{Step 100 ($P=1,963$)}
	\end{subfigure}
	\begin{subfigure}{1.0\textwidth}
		\centering
		\includegraphics[scale=0.3]{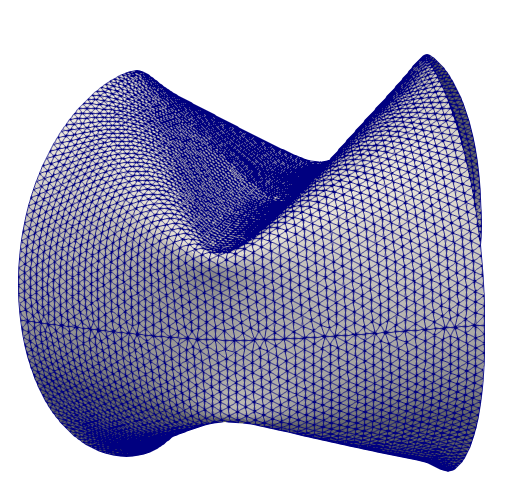}
		\caption{Step 200 ($P=5,620$)}
	\end{subfigure}
	\begin{subfigure}{1.0\textwidth}
		\centering
		\includegraphics[scale=0.3]{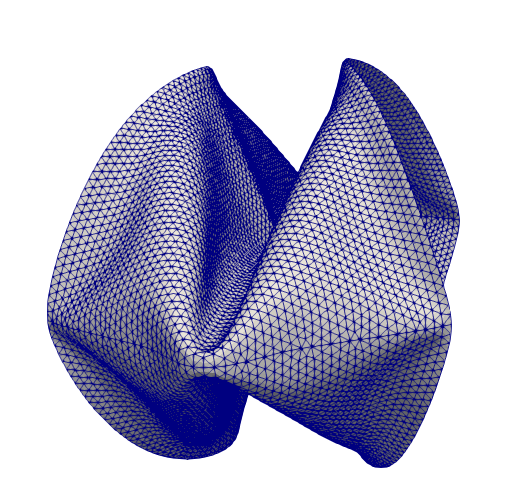}
		\caption{Step 368 (last step, $P=48,529$)}
	\end{subfigure}
	\caption{Deformed shapes of the pinched cylindrical shell at three stages of the analysis. 
	\label{fig:pinched_shell_steps}}
\end{figure}

Fig. \ref{fig:pinched_cylinder_comparison} presents the
displacements at observed points A and B of the shell in comparison with the
reference \cite{sze_popular_2004}. 
Similar to the previous examples,  second order
approximation is employed in the through-thickness direction and various orders
are investigated for upper and lower faces. As can be seen, there is good agreement with the reference results. Meanwhile, the
shapes of the pinched shell at three stages of the analysis are shown in Fig.
\ref{fig:pinched_shell_steps}, where the large and complex deformation at the
middle of the cylinder (point B) is captured. The analysis ends when the two sides of the cylinder meet.

\section{Concluding remarks}
\label{sec:conclusion}

In this study, the formulation of a prism element using shape 
functions that are hierarchical, heterogeneous, and anisotropic is presented 
for solving complex nonlinear behaviour of solve solid shell problems.
Unlike classical degenerated shell theory, the solid shell theory requires only 
displacement DOFs, enabling them to be combined with classical solid elements if
required and avoids complex description of boundary conditions. Furthermore,  
stretching in the through-thickness direction is captured and general 3D 
constitutive equations can be used without simplification or modification.

Additionally, element of this type enables modelling of coupled multi-physics
problems as well as utilisation of more sophisticated hyperelastic models since
the material is defined in a standard 3D setting.

The hierarchical properties of the shape functions enable
an efficient iterative solver to be employed, that is tailored for \emph{hp}-adaptive 
implementation. The heterogeneous property of the shape functions enables the 
use of arbitrary orders of approximation that can be set independently for each 
mesh entity, i.e. edge, triangle, quadrilateral, and prism. The anisotropic 
property permits the approximation functions on the upper and lower faces
to be defined independently from approximation through the thickness
of the element. Thus, locking problems, that are 
common in solid shell formulations, are alleviated simply by increasing the order 
of approximation, without resorting to enhancements, such as reduced integration 
or additional assumed natural strain or enhanced strain fields. 

Looking forward, through layering of the the proposed prism elements
can be easily applied tor modelling laminated structures and the delamination
is captured by adding standard cohesive/interface elements between layers. 
It is also straightforward to add 
dissipative processes such as plasticity or fracture.

The formulation has been implemented in the \code{MoFEM} library \cite{kaczmarczyk_mofem_2020}. Both the code and the data
for the numerical examples are open-source \cite{mofem_solid_shell_code_and_data}.

\section*{Acknowledgements}

The authors gratefully acknowledge the support provided by the EPSRC Strategic
Support Package: Engineering of Active Materials by Multiscale/Multiphysics
Computational Mechanics - EP/R008531/1.

\bibliographystyle{elsarticle-num}
\bibliography{solid_shell_manuscript_pre_print}

\end{document}